# THE DANTZIG SELECTOR: STATISTICAL ESTIMATION WHEN $p$ IS MUCH LARGER THAN $n$[1]

By Emmanuel Candes[2] and Terence Tao[3]

*California Institute of Technology and University of California, Los Angeles*

In many important statistical applications, the number of variables or parameters $p$ is much larger than the number of observations $n$. Suppose then that we have observations $y = X\beta + z$, where $\beta \in \mathbf{R}^p$ is a parameter vector of interest, $X$ is a data matrix with possibly far fewer rows than columns, $n \ll p$, and the $z_i$'s are i.i.d. $N(0, \sigma^2)$. Is it possible to estimate $\beta$ reliably based on the noisy data $y$?

To estimate $\beta$, we introduce a new estimator—we call it the *Dantzig selector*—which is a solution to the $\ell_1$-regularization problem

$$\min_{\tilde{\beta} \in \mathbf{R}^p} \|\tilde{\beta}\|_{\ell_1} \quad \text{subject to} \quad \|X^* r\|_{\ell_\infty} \leq (1 + t^{-1})\sqrt{2 \log p} \cdot \sigma,$$

where $r$ is the residual vector $y - X\tilde{\beta}$ and $t$ is a positive scalar. We show that if $X$ obeys a uniform uncertainty principle (with unit-normed columns) and if the true parameter vector $\beta$ is sufficiently sparse (which here roughly guarantees that the model is identifiable), then with very large probability,

$$\|\hat{\beta} - \beta\|_{\ell_2}^2 \leq C^2 \cdot 2 \log p \cdot \left( \sigma^2 + \sum_i \min(\beta_i^2, \sigma^2) \right).$$

Our results are nonasymptotic and we give values for the constant $C$. Even though $n$ may be much smaller than $p$, our estimator achieves a loss within a logarithmic factor of the ideal mean squared error one

Received August 2005; revised March 2006.

[1]Discussed in 10.1214/009053607000000424, 10.1214/009053607000000433, 10.1214/009053607000000442, 10.1214/009053607000000451, 10.1214/009053607000000460 and 10.1214/009053607000000479; rejoinder at 10.1214/009053607000000532.

[2]Supported in part by NSF Grant DMS-01-40698 and by an Alfred P. Sloan Fellowship.

[3]Supported in part by a grant from the Packard Foundation.

would achieve with an *oracle* which would supply perfect information about which coordinates are nonzero, and which were above the noise level.

In multivariate regression and from a model selection viewpoint, our result says that it is possible nearly to select the best subset of variables by solving a very simple convex program, which, in fact, can easily be recast as a convenient linear program (LP).

## 1. Introduction.

In many important statistical applications, the number of variables or parameters $p$ is now much larger than the number of observations $n$. In radiology and biomedical imaging, for instance, one is typically able to collect far fewer measurements about an image of interest than the unknown number of pixels. Examples in functional MRI and tomography all come to mind. High dimensional data frequently arise in genomics. Gene expression studies are a typical example: a relatively low number of observations (in the tens) is available, while the total number of genes assayed (and considered as possible regressors) is easily in the thousands. Other examples in statistical signal processing and nonparametric estimation include the recovery of a continuous-time curve or surface from a finite number of noisy samples. Estimation in this setting is generally acknowledged as an important challenge in contemporary statistics; see the recent conference held in Leiden, The Netherlands (September 2002), "On high-dimensional data $p \gg n$ in mathematical statistics and bio-medical applications." It is believed that progress may have the potential for impact across many areas of statistics [30].

In many research fields then, scientists work with data matrices with many variables $p$ and comparably few observations $n$. This paper is about this important situation, and considers the problem of estimating a parameter $\beta \in \mathbf{R}^p$ from the linear model

$$(1.1) \qquad\qquad y = X\beta + z;$$

$y \in \mathbf{R}^n$ is a vector of observations, $X$ is an $n \times p$ predictor matrix, and $z$ a vector of stochastic measurement errors. Unless specified otherwise, we will assume that $z \sim N(0, \sigma^2 I_n)$ is a vector of independent normal random variables, although it is clear that our methods and results may be extended to other distributions. Throughout this paper, we will of course typically assume that $p$ is much larger than $n$.

### 1.1. *Uniform uncertainty principles and the noiseless case.*

At first, reliably estimating $\beta$ from $y$ may seem impossible. Even in the noiseless case, one may wonder how one could possibly do this, as one would need to solve an underdetermined system of equations with fewer equations than unknowns. But suppose now that $\beta$ is known to be structured in the sense that it is sparse or compressible. For example, suppose that $\beta$ is $S$-sparse so that only $S$ of its entries are nonzero. This premise radically changes the problem, making the search for solutions feasible. In fact, [13] showed that in the noiseless case, one could actually recover $\beta$ *exactly* by solving the convex program ($\|\tilde{\beta}\|_{\ell_1} := \sum_{i=1}^{p} |\tilde{\beta}_i|$)

$$(1.2) \qquad (P_1) \qquad \min_{\tilde{\beta} \in \mathbf{R}^p} \|\tilde{\beta}\|_{\ell_1} \quad \text{subject to} \quad X\tilde{\beta} = y,$$

provided that the matrix $X \in \mathbf{R}^{n \times p}$ obeys a *uniform uncertainty principle*. (The program $(P_1)$ can even be recast as a linear program.) That is, $\ell_1$-minimization finds without error both the location and amplitudes—which we emphasize are a priori completely unknown—of the nonzero components of the vector $\beta \in \mathbf{R}^p$. We also refer the reader to [20, 24] and [27] for inspiring early results.

To understand the exact recovery phenomenon, we introduce the notion of uniform uncertainty principle (UUP) proposed in [14] and refined in [13]. This principle will play an important role throughout, although we emphasize that this paper *is not* about the exact recovery of noiseless data. The UUP essentially states that the $n \times p$ measurement or design matrix $X$ obeys a "restricted isometry hypothesis." Let $X_T$, $T \subset \{1, \ldots, p\}$, be the $n \times |T|$ submatrix obtained by extracting the columns of $X$ corresponding to the indices in $T$; then [13] defines the $S$-restricted isometry constant $\delta_S$ of $X$ which is the smallest quantity such that

$$(1.3) \qquad (1 - \delta_S)\|c\|_{\ell_2}^2 \le \|X_T c\|_{\ell_2}^2 \le (1 + \delta_S)\|c\|_{\ell_2}^2$$

for all subsets $T$ with $|T| \le S$ and coefficient sequences $(c_j)_{j \in T}$. This property essentially requires that every set of columns with cardinality less than $S$ approximately behaves like an orthonormal system. It was shown (also in [13]) that if $S$ obeys

$$(1.4) \qquad \delta_S + \delta_{2S} + \delta_{3S} < 1,$$

then solving $(P_1)$ recovers *any* sparse signal $\beta$ with support size obeying $|T| \le S$.

Actually, [13] derived a slightly stronger result. Introduce the $S, S'$-*restricted orthogonality constants* $\theta_{S,S'}$ for $S + S' \le p$ to be the smallest quantities such that

$$(1.5) \qquad |\langle X_T c, X_{T'} c' \rangle| \le \theta_{S,S'} \cdot \|c\|_{\ell_2} \|c'\|_{\ell_2}$$

holds for all *disjoint* sets $T, T' \subseteq \{1, \ldots, p\}$ of cardinality $|T| \le S$ and $|T'| \le S'$. Small values of restricted orthogonality constants indicate that disjoint subsets of covariates span nearly orthogonal subspaces. Then the authors showed that the recovery is exact, provided

$$(1.6) \qquad \delta_S + \theta_{S,S} + \theta_{S,2S} < 1,$$

which is a little better since it is not hard to see that $\delta_{S+S'} - \delta_{S'} \le \theta_{S,S'} \le \delta_{S+S'}$ for $S' \ge S$ ([13], Lemma 1.1).

1.2. *Uniform uncertainty principles and statistical estimation.* Any real-world sensor or measurement device is subject to at least a small amount of noise. And now one asks whether it is possible to reliably estimate the parameter $\beta \in \mathbf{R}^p$ from the noisy data $y \in \mathbf{R}^n$ and the model (1.1). Frankly, this may seem like an impossible task. How can one hope to estimate $\beta$, when, in addition to having too few observations, these are also contaminated with noise?

To estimate $\beta$ with noisy data, we consider, nevertheless, solving the convex program

$$(1.7) \quad (DS) \quad \min_{\tilde{\beta} \in \mathbf{R}^p} \|\tilde{\beta}\|_{\ell_1} \quad \text{subject to} \quad \|X^*r\|_{\ell_\infty} := \sup_{1 \le i \le p} |(X^*r)_i| \le \lambda_p \cdot \sigma$$

for some $\lambda_p > 0$, where $r$ is the vector of residuals

$$(1.8) \qquad\qquad\qquad r = y - X\tilde{\beta}.$$

In other words, we seek an estimator $\hat{\beta}$ with minimum complexity (as measured by the $\ell_1$-norm) among all objects that are consistent with the data. The constraint on the residual vector imposes that for each $i \in \{1, \ldots, p\}$, $|(X^*r)_i| \le \lambda_p \cdot \sigma$, and guarantees that the residuals are within the noise level. As we shall see later, this proposal makes sense provided that the columns of $X$ have the same Euclidean size and in this paper we will always assume they are unit-normed; our results would equally apply to matrices with different column sizes—one would need to change the right-hand side to $|(X^*r)_i|$ less or equal to $\lambda_p \cdot \sigma$ times the Euclidean norm of the $i$th column of $X$, or to $|(X^*r)_i| \le \sqrt{1 + \delta_1} \cdot \lambda_p \cdot \sigma$ since all the columns have norm less than $\sqrt{1 + \delta_1}$.

There are many reasons why one would want to constrain the size of the correlated residual vector $X^*r$ rather than the size of the residual vector $r$. Suppose that an orthonormal transformation is applied to the data, giving $y' = Uy$, where $U^*U$ is the identity. Clearly, a good estimation procedure for estimating $\beta$ should not depend upon $U$ (after all, one could apply $U^*$ to return to the original problem). It turns out that the estimation procedure (1.7) is actually invariant with respect to orthonormal transformations applied to the data vector since the feasible region is invariant: $(UX)^T(UX\tilde{\beta} - Uy) = X^*(X\tilde{\beta} - y)$. In contrast, had we defined the feasibility region with $\sup_i |r_i|$ being smaller than a fixed threshold, then the estimation procedure would not be invariant. There are other reasons aside from this. One of them is that we would obviously want to include in the model explanatory variables that are highly correlated with the data $y$. Consider the situation in which a residual vector is equal to a column $X^i$ of the design matrix $X$. Suppose, for simplicity, that the components of $X^i$ all have about the same size, that is, about $1/\sqrt{n}$, and assume that $\sigma$ is slightly larger than $1/\sqrt{n}$. Had we used a constraint of the form $\sup_i |r_i| \le \lambda_n \sigma$ (with perhaps $\lambda_n$ of size about $\sqrt{2 \log n}$), the vector of residuals would be feasible, which does not make any sense. In contrast, such a residual vector would not be feasible for (1.7) for reasonable values of the noise level, and the $i$th variable would be rightly included in the model.

Again, the program $(DS)$ is convex, and can easily be recast as a linear program (LP),

$$\min \sum_i u_i \quad \text{subject to} \quad -u \le \tilde{\beta} \le u \quad \text{and}$$

(1.9)

$$-\lambda_p \sigma \mathbf{1} \le X^*(y - X\tilde{\beta}) \le \lambda_p \sigma \mathbf{1},$$

where the optimization variables are $u, \tilde{\beta} \in R^p$, and $\mathbf{1}$ is a $p$-dimensional vector of ones. Hence, our estimation procedure is computationally tractable; see Section 4.4 for details. There is indeed a growing family of ever more efficient algorithms for solving such problems (even for problems with tens or even hundreds of thousands of observations) [8].

We call the estimator (1.7) the *Dantzig selector*; with this name, we intend to pay tribute to the father of linear programming who passed away while we were finalizing this manuscript, and to underscore that the convex program $(DS)$ is effectively a variable selection technique.

The first result of this paper is that the Dantzig selector is surprisingly accurate.

THEOREM 1.1. *Suppose $\beta \in \mathbf{R}^p$ is any $S$-sparse vector of parameters obeying $\delta_{2S} + \theta_{S,2S} < 1$. Choose $\lambda_p = \sqrt{2 \log p}$ in (1.7). Then with large probability, $\hat{\beta}$ obeys*

(1.10)
$$\|\hat{\beta} - \beta\|_{\ell_2}^2 \le C_1^2 \cdot (2 \log p) \cdot S \cdot \sigma^2,$$

*with $C_1 = 4/(1 - \delta_S - \theta_{S,2S})$. Hence, for small values of $\delta_S + \theta_{S,2S}$, $C_1 \approx 4$. For concreteness, if one chooses $\lambda_p := \sqrt{2(1+a) \log p}$ for each $a \ge 0$, the bound holds with probability exceeding $1 - (\sqrt{\pi \log p} \cdot p^a)^{-1}$ with the proviso that $\lambda_p^2$ substitutes $2 \log p$ in (1.10).*

We will discuss the condition $\delta_{2S} + \theta_{S,2S} < 1$ later but, for the moment, observe that (1.10) describes a striking phenomenon: not only are we able to reliably estimate the vector of parameters from limited observations, but the mean squared error is simply proportional—up to a logarithmic factor—to the true number of unknowns times the noise level $\sigma^2$. What is of interest here is that one can achieve this feat by solving a simple linear program. Moreover, and ignoring the log-like factor, statistical common sense tells us that (1.10) is, in general, unimprovable.

To see why this is true, suppose one had available an oracle letting us know in advance the location of the $S$ nonzero entries of the parameter vector, that is, $T_0 := \{i : \beta_i \neq 0\}$. That is, in the language of model selection, one would know the right model ahead of time. Then one could use this valuable information and construct an ideal estimator $\beta^\star$ by using the least-squares projection

$$\beta_{T_0}^\star = (X_{T_0}^T X_{T_0})^{-1} X_{T_0}^T y,$$

where $\beta_{T_0}^\star$ is the restriction of $\beta^\star$ to the set $T_0$, and set $\beta^\star$ to zero outside of $T_0$. (At times, we will abuse notation and also let $\beta_I$ be the truncated vector equal to $\beta_i$ for $i \in I$ and zero otherwise.) Clearly,

$$\beta^\star = \beta + (X_{T_0}^T X_{T_0})^{-1} X_{T_0}^T z$$

and

$$\mathbf{E}\|\beta^\star - \beta\|_{\ell_2}^2 = \mathbf{E}\|(X_{T_0}^T X_{T_0})^{-1} X_{T_0}^T z\|_{\ell_2}^2 = \sigma^2 \operatorname{Tr}((X_{T_0}^T X_{T_0})^{-1}).$$

Now since all the eigenvalues of $X_{T_0}^T X_{T_0}$ belong to the interval $[1 - \delta_S, 1 + \delta_S]$, the ideal expected mean squared error would obey

$$\mathbf{E}\|\beta^\star - \beta\|_{\ell_2}^2 \geq \frac{1}{1 + \delta_S} \cdot S \cdot \sigma^2.$$

Hence, Theorem 1.1 says that the minimum $\ell_1$ estimator achieves a loss within a logarithmic factor of the ideal mean squared error; the logarithmic factor is the price we pay for adaptivity, that is, for not knowing ahead of time where the nonzero parameter values actually are.

In short, the recovery procedure, although extremely nonlinear, is *stable in the presence of noise*. This is especially interesting because the matrix $X$ in (1.1) is rectangular; it has many more columns than rows. As such, most of its singular values are zero. In solving $(DS)$, we are essentially trying to invert the action of $X$ on our hidden $\beta$ in the presence of noise. The fact that this matrix inversion process keeps the perturbation from "blowing up"—even though it is severely ill-posed—is perhaps unexpected.

Presumably, our result would be especially interesting if one could estimate the order of $n$ parameters with as few as $n$ observations. That is, we would like the condition $\delta_{2S} + \theta_{S,2S} < 1$ to hold for very large values of $S$, for example, as close as possible to $n$ (note that for $2S > n$, $\delta_{2S} \geq 1$ since any submatrix with more than $n$ columns must be singular, which implies that in any event $S$ must be less than $n/2$). Now, this paper is part of a larger body of work [11, 13, 14] which shows that, for "generic" or random design matrices $X$, the condition holds for very significant values of $S$. Suppose, for instance, that $X$ is a random matrix with i.i.d. Gaussian entries. Then with overwhelming probability, the condition holds for $S = O(n/\log(p/n))$. In other words, this setup only requires $O(\log(p/n))$ observations per nonzero parameter value; for example, when $n$ is a nonnegligible fraction of $p$, one only needs a handful of observations per nonzero coefficient. In practice, this number is quite small, as few as 5 or 6 observations per unknown generally suffice (over a large range of the ratio $p/n$); see Section 4. Many design matrices have a similar behavior and Section 2 discusses a few of these.

As an aside, it is interesting to note that, for $S$ obeying the condition of the theorem, the reconstruction from noiseless data ($\sigma = 0$) is exact and that our condition is slightly better than (1.6).

1.3. *Oracle inequalities.* Theorem 1.1 is certainly noticeable but there are instances, however, in which it may still be a little naive. Suppose, for example, that $\beta$ is very small so that $\beta$ is well below the noise level, that is, $|\beta_i| \ll \sigma$ for all $i$. Then with this information we could set $\hat{\beta} = 0$, and the squared error loss would then simply be $\sum_{i=1}^{p} |\beta_i|^2$, which may potentially be much smaller than $\sigma^2$ times the number of nonzero coordinates of $\beta$. In some sense, this is a situation in which the squared bias is much smaller than the variance.

A more ambitious proposal might then ask for a near-optimal trade-off *coordinate by coordinate*. To explain this idea, suppose, for simplicity, that $X$ is the identity matrix so that $y \sim N(\beta, \sigma^2 I_p)$. Suppose then that we had available an *oracle* letting us know ahead of time which coordinates of $\beta$ are significant, that is, the set of indices for which $|\beta_i| > \sigma$. Then equipped with this oracle, we would set $\beta_i^\star = y_i$ for each index in the significant set and $\beta_i^\star = 0$ otherwise. The expected mean squared error of this ideal estimator is then

$$(1.11) \qquad \mathbf{E}\|\beta^\star - \beta\|_{\ell_2}^2 = \sum_{i=1}^{p} \min(\beta_i^2, \sigma^2).$$

Here and below, we will refer to (1.11) as the ideal MSE. As is well known, thresholding rules with threshold level at about $\sqrt{2\log p} \cdot \sigma$ achieve the ideal MSE to within a multiplicative factor proportional to $\log p$ [21, 22].

In the context of the linear model, we might think about the ideal estimation as follows: consider the least-squares estimator $\hat{\beta}_I = (X_I^T X_I)^{-1} X_I^T y$ as before and consider the ideal least-squares estimator $\beta^\star$ which minimizes the expected mean squared error

$$\beta^\star = \underset{I \subset \{1, \dots, p\}}{\arg\min} \mathbf{E}\|\beta - \hat{\beta}_I\|_{\ell_2}^2.$$

In other words, one would fit all least-squares models and rely on an oracle to tell us which model to choose. This is *ideal* because we can of course not evaluate $\mathbf{E}\|\beta - \hat{\beta}_I\|_{\ell_2}^2$ since we do not know $\beta$ (we are trying to estimate it after all). But we can view this as a benchmark and ask whether any real estimator would obey

$$(1.12) \qquad \|\hat{\beta} - \beta\|_{\ell_2}^2 = O(\log p) \cdot \mathbf{E}\|\beta - \beta^\star\|_{\ell_2}^2$$

with large probability.

In some sense, (1.11) is a proxy for the ideal risk $\mathbf{E}\|\beta - \beta^\star\|_{\ell_2}^2$. Indeed, let $I$ be a fixed subset of indices and consider regressing $y$ onto this subset (we again denote by $\beta_I$ the restriction of $\beta$ to the set $I$). The error of this estimator is given by

$$\|\hat{\beta}_I - \beta\|_{\ell_2}^2 = \|\hat{\beta}_I - \beta_I\|_{\ell_2}^2 + \|\beta_I - \beta\|_{\ell_2}^2.$$

The first term is equal to

$$\hat{\beta}_I - \beta_I = (X_I^T X_I)^{-1} X_I^T X \beta_{I^c} + (X_I^T X_I)^{-1} X_I^T z,$$

and its expected mean squared error is given by the formula

$$\mathbf{E}\|\hat{\beta}_I - \beta_I\|^2 = \|(X_I^T X_I)^{-1} X_I^T X \beta_{I^c}\|_{\ell_2}^2 + \sigma^2 \operatorname{Tr}((X_I^T X_I)^{-1}).$$

Thus, this term obeys

$$\mathbf{E}\|\hat{\beta}_I - \beta_I\|^2 \geq \frac{1}{1 + \delta_{|I|}} \cdot |I| \cdot \sigma^2$$

for the same reasons as before. In short, for all sets $I$, $|I| \leq S$, with $\delta_S < 1$, say,

$$\mathbf{E}\|\hat{\beta}_I - \beta\|^2 \geq \frac{1}{2} \cdot \left( \sum_{i \in I^c} \beta_i^2 + |I| \cdot \sigma^2 \right),$$

which gives that the ideal mean squared error is bounded below by

$$\mathbf{E}\|\beta^\star - \beta\|_{\ell_2}^2 \geq \frac{1}{2} \cdot \min_I \left( \sum_{i \in I^c} \beta_i^2 + |I| \cdot \sigma^2 \right) = \frac{1}{2} \cdot \sum_i \min(\beta_i^2, \sigma^2).$$

In that sense, the ideal risk is lower bounded by the proxy (1.12). As we have seen, the proxy is meaningful since it has a natural interpretation in terms of the ideal squared bias and variance,

$$\sum_i \min(\beta_i^2, \sigma^2) = \min_{I \subset \{1,\ldots,p\}} \|\beta - \beta_I\|_{\ell_2}^2 + |I| \cdot \sigma^2.$$

This raises a fundamental question: given data $y$ and the linear model (1.1), not knowing anything about the significant coordinates of $\beta$ and not being able to observe directly the parameter values, can we design an estimator which nearly achieves (1.12)? Our main result is that the Dantzig selector (1.7) does just that.

THEOREM 1.2. *Choose* $t > 0$ *and set* $\lambda_p := (1 + t^{-1})\sqrt{2 \log p}$ *in* (1.7). *Then if* $\beta$ *is* $S$-sparse with $\delta_{2S} + \theta_{S,2S} < 1 - t$, *our estimator obeys*

$$(1.13) \qquad \|\hat{\beta} - \beta\|_{\ell_2}^2 \leq C_2^2 \cdot \lambda_p^2 \cdot \left( \sigma^2 + \sum_{i=1}^{p} \min(\beta_i^2, \sigma^2) \right)$$

*with large probability [the probability is as before for $\lambda_p := (\sqrt{1+a} + t^{-1}) \times \sqrt{2\log p}$]. Here, $C_2$ may only depend on $\delta_{2S}$ and $\theta_{S,2S}$; see below.*

We emphasize that (1.13) is nonasymptotic and our analysis actually yields explicit constants. For instance, we also prove that

$$C_2 = 2\frac{C_0}{1-\delta-\theta} + 2\frac{\theta(1+\delta)}{(1-\delta-\theta)^2} + \frac{1+\delta}{1-\delta-\theta}$$

and

$$(1.14) \qquad C_0 := 2\sqrt{2}\left(1 + \frac{1-\delta^2}{1-\delta-\theta}\right) + (1+1/\sqrt{2})\frac{(1+\delta)^2}{1-\delta-\theta},$$

where above and below, we put $\delta := \delta_{2S}$ and $\theta := \theta_{S,2S}$ for convenience. For $\delta$ and $\theta$ small, $C_2$ is close to

$$C_2 \approx 2(4\sqrt{2} + 1 + 1/\sqrt{2}) + 1 \le 16.$$

The condition imposing $\delta_{2S} + \theta_{S,2S} < 1$ (or less than $1-t$) has a rather natural interpretation in terms of model identifiability. Consider a rank deficient submatrix $X_{T\cup T'}$ with $2S$ columns (lowest eigenvalue is $0 = 1 - \delta_{2S}$), and with indices in $T$ and $T'$, each of size $S$. Then there is a vector $h$ obeying $Xh = 0$ and which can be decomposed as $h = \beta - \beta'$, where $\beta$ is supported on $T$ and likewise for $\beta'$; that is,

$$X\beta = X\beta'.$$

In short, this says that the model is not identifiable since both $\beta$ and $\beta'$ are $S$-sparse. In other words, we need $\delta_{2S} < 1$. The requirement $\delta_{2S} + \theta_{S,2S} < 1$ (or less than $1-t$) is only slightly stronger than the identifiability condition, roughly two times stronger. It puts a lower bound on the singular values of submatrices and, in effect, prevents situations where multicollinearity between competitive subsets of predictors could occur.

1.4. *Ideal model selection by linear programming.* Our estimation procedure is of course an implicit method for choosing a desirable subset of predictors, based on the noisy data $y = X\beta + z$, from among all subsets of variables. As the reader will see, there is nothing in our arguments that requires $p$ to be larger than $n$ and, thus, the Dantzig selector can be construed as a very general variable selection strategy—hence, the name.

There is of course a huge literature on model selection, and many procedures motivated by a wide array of criteria have been proposed over the years—among which [1, 7, 26, 31, 36]. By and large, the most commonly discussed approach—the "canonical selection procedure" according to [26]—is defined as

$$(1.15) \qquad \underset{\tilde{\beta} \in \mathbf{R}^p}{\arg\min} \|y - X\tilde{\beta}\|_{\ell_2}^2 + \Lambda \cdot \sigma^2 \cdot \|\tilde{\beta}\|_{\ell_0}, \qquad \|\tilde{\beta}\|_{\ell_0} := |\{i : \tilde{\beta}_i \ne 0\}|,$$

which best trades-off between the goodness of fit and the complexity of the model, the so-called bias and variance terms. Popular selection procedures such as AIC, $C_p$, BIC and RIC are all of this form with different values of the parameter $\Lambda$; see also [2, 4, 5, 6, 7] for related proposals. To make a long story short, model selection is an important area in part because of the thousands of people routinely fitting large linear models or designing statistical experiments. As such, it has and still receives a lot of attention, and progress in this field is likely to have a large impact. Now despite the size of the current literature, we believe there are two critical problems in this field:

- First, finding the minimum of (1.15) is in general *NP*-hard [32]. To the best of our knowledge, solving this problem essentially requires exhaustive searches over all subsets of columns of $X$, a procedure which clearly is combinatorial in nature and has exponential complexity since for $p$ of size about $n$, there are about $2^p$ such subsets. (We are of course aware that in a few exceptional circumstances, e.g., when $X$ is an orthonormal matrix, the solution is computationally feasible and given by thresholding rules [7, 21].)

    In other words, solving the model selection problem might be possible only when $p$ ranges in the few dozens. This is especially problematic when one considers that we now live in a "data-driven" era marked by ever larger datasets.

- Second, estimating $\beta$ and $X\beta$—especially when $p$ is larger than $n$—are two very different problems. Whereas there is an extensive literature about the problem of estimating $X\beta$, the quantitative literature about the equally important problem of estimating $\beta$ in the modern setup where $p$ is not small compared to $n$ is scarce; see [25]. For completeness, important and beautiful results about the former problem (estimating $X\beta$) include the papers [3, 4, 6, 7, 23, 26].

In recent years, researchers have of course developed alternatives to overcome these computational difficulties, and we would like to single out the popular lasso also known as Basis Pursuit [15, 38], which relaxes the counting norm $\|\tilde{\beta}\|_{\ell_0}$ into the convex $\ell_1$-norm $\|\tilde{\beta}\|_{\ell_1}$. Notwithstanding the novel and exciting work of [28] on the persistence of the lasso for variable selection in high dimensions—which again is about estimating $X\beta$ and not $\beta$—not much is yet known about the performance of such strategies although they seem to work well in practice; for example, see [35].

Against this background, our work clearly marks a significant departure from the current literature, both in terms of what it achieves and of its methods. Indeed, our paper introduces a method for selecting variables based on linear programming, and obtains decisive quantitative results in fairly general settings.

1.5. *Extension to nearly sparse parameters.* We have considered thus far the estimation of sparse parameter vectors, that is, with a number $S$ of nonzero entries obeying $\delta_{2S} + \theta_{S,2S}$. We already explained that this condition is in some sense necessary as otherwise one might have an "aliasing" problem, a situation in which $X\beta \approx X\beta'$, although $\beta$ and $\beta'$ might be completely different. However, extensions of our results to nonsparse objects are possible provided that one imposes other types of constraints to remove the possibility of strong aliasing.

Many such constraints may exist and we consider one of them which imposes some decay condition on the entries of $\beta$. Rearrange the entries of $\beta$ by decreasing order of magnitude $|\beta_{(1)}| \geq |\beta_{(2)}| \geq \cdots \geq |\beta_{(p)}|$ and suppose the $k$th largest entry obeys

$$(1.16) \qquad |\beta_{(k)}| \leq R \cdot k^{-1/s},$$

for some positive $R$ and $s \leq 1$, say. Can we show that our estimator achieves an error close to the proxy (1.11)? The first observation is that, to mimic this proxy, we need to be able to estimate reliably all the coordinates which are significantly above the noise level, that is, roughly such that $|\beta_i| \geq \sigma$. Let $S = |\{i : |\beta_i| > \sigma\}|$. Then if $\delta_{2S} + \theta_{S,2S} < 1$, this might be possible, but otherwise, we may simply not have enough observations to estimate that many coefficients. The second observation is that for $\beta \in \mathbf{R}^p$ obeying (1.16),

$$(1.17) \qquad \sum_i \min(\beta_i^2, \sigma^2) = S \cdot \sigma^2 + \sum_{i \geq S+1} |\beta_{(i)}|^2 \leq C \cdot (S \cdot \sigma^2 + R^2 S^{-2r})$$

with $r = 1/s - 1/2$. With this in mind, we have the following result.

THEOREM 1.3. *Suppose $\beta \in \mathbf{R}^p$ obeys* (1.16) *and let $S_*$ be fixed such that $\delta_{2S_*} + \theta_{S_*,2S_*} < 1$. Choose $\lambda_p$ as in Theorem 1.1. Then $\hat{\beta}$ obeys*

$$(1.18) \qquad \|\hat{\beta} - \beta\|_{\ell_2}^2 \leq \min_{1 \leq S \leq S_*} C_3 \cdot 2 \log p \cdot (S \cdot \sigma^2 + R^2 S^{-2r})$$

*with large probability.*

Note that for each $\beta$ obeying (1.16), $|\{i : |\beta_i| > \sigma\}| \leq (R/\sigma)^{1/s}$. Then if $S_* \geq (R/\sigma)^{1/s}$, it is not hard to see that (1.18) becomes

$$(1.19) \qquad \|\hat{\beta} - \beta\|_{\ell_2}^2 \leq O(\log p) \cdot R^{2/(2r+1)} \cdot (\sigma^2)^{2r/(2r+1)},$$

which is the well-known minimax rate for classes of objects exhibiting the decay (1.16). Even though we have $n \ll p$, the Dantzig selector recovers the minimax rate that one would get if we were able to measure *all* the coordinates of $\beta$ *directly* via $\tilde{y} \sim N(\beta, \sigma^2 I_p)$. In the case where $S_* \leq (R/\sigma)^{1/s}$, the method saturates because we do not have enough data to recover the minimax rate, and can only guarantee a squared loss of about $O(\log p)(R^2 S_*^{-2r} + S_* \cdot \sigma^2)$. Note, however, that the error is well controlled.

1.6. *Variations and other extensions.* When $X$ is an orthogonal matrix, the Dantzig selector $\hat{\beta}$ is then the $\ell_1$-minimizer subject to the constraint $\|X^*y - \hat{\beta}\|_{\ell_\infty} \le \lambda_p \cdot \sigma$. This implies that $\hat{\beta}$ is simply the soft-thresholded version of $X^*y$ at level $\lambda_p \cdot \sigma$; thus,

$$\hat{\beta}_i = \max(|(X^*y)_i| - \lambda_p \cdot \sigma, 0)\operatorname{sgn}((X^*y)_i).$$

In other words, $X^*y$ is shifted toward the origin. In Section 4 we will see that for arbitrary $X$'s the method continues to exhibit a soft-thresholding type of behavior and as a result, may slightly underestimate the true value of the nonzero parameters.

There are several simple methods which can correct for this bias and increase performance in practical settings. We consider one of these based on a two-stage procedure:

1. Estimate $I = \{i : \beta_i \ne 0\}$ with $\hat{I} = \{i : \hat{\beta}_i \ne 0\}$ with $\beta$ as in (1.7) (or, more generally, with $\hat{I} = \{i : |\hat{\beta}_i| > \alpha \cdot \sigma\}$ for some $\alpha \ge 0$).
2. Construct the estimator

$$(1.20) \qquad\qquad \hat{\beta}_{\hat{I}} = (X_{\hat{I}}^T X_{\hat{I}})^{-1} X_{\hat{I}}^T y,$$

   and set the other coordinates to zero.

Hence, we rely on the Dantzig selector to estimate the model $I$, and construct a new estimator by regressing the data $y$ onto the model $\hat{I}$. We will refer to this variation as the Gauss–Dantzig selector. As we will see in Section 4, this recenters the estimate and generally yields higher statistical accuracy. We anticipate that all our theorems hold with some such variations.

Although we prove our main results in the case where $z$ is a vector of i.i.d. Gaussian variables, our methods and results would certainly extend to other noise distributions. The key is to constrain the residuals so that the true vector $\beta$ is feasible for the optimization problem. In details, this means that we need to set $\lambda_p$ so that $Z^* = \sup_i |\langle X^i, z \rangle|$ is less than $\lambda_p \sigma$ with large probability. When $z \sim N(0, \sigma^2 I_n)$, this is achieved for $\lambda_p = \sqrt{2 \log p}$, but one could compute other thresholds for other types of zero-mean distributions and derive results similar to those introduced in this paper. In general setups, one would perhaps want to be more flexible and have thresholds depending upon the column index. For example, one could declare that $r$ is feasible if $\sup_i |\langle X^i, r \rangle| / \lambda_p^i$ is below some fixed threshold.

1.7. *Organization of the paper.* The paper is organized as follows. We begin by discussing the implications of this work for experimental design in Section 2. We prove our main results, namely, Theorems 1.1, 1.2 and 1.3, in Section 3. Section 4 introduces numerical experiments showing that our approach is effective in practical applications. Finally, Section 5 closes the paper with a short summary of our findings and of their consequences for model selection, and with a discussion of other related work in Section 5.2. Finally, the Appendix provides proofs of key lemmas supporting the proof of Theorem 1.2.

**2. Significance for experimental design.** Before we begin proving our main results, we would like to explain why our method might be of interest to anyone seeking to measure or sense a sparse high-dimensional vector using as few measurements as possible. In the noiseless case, our earlier results showed that if $\beta$ is $S$-sparse, then it can be reconstructed *exactly* from $n$ measurements $y = X\beta$, provided that $\delta + \theta < 1$ [11, 13]. These were later extended to include wider classes of objects, that is, the so called *compressible* objects. Against this background, our results show that the Dantzig selector is robust against measurement errors (no realistic measuring device can provide infinite precision), thereby making it well suited for practical applications.

2.1. *Random matrices and designs.* An interesting aspect of this theory is that random matrices $X$ are in some sense ideal for recovering an object from a few projections. For example, if $X$ is a properly normalized Gaussian matrix with i.i.d. entries, then the conditions of our theorems hold with

$$(2.1) \qquad S \asymp n/\log(p/n)$$

with overwhelming probability [13, 14, 18, 37]. The same relation is also conjectured to be true for other types of random matrices such as normalized binary arrays with i.i.d. entries taking values $\pm 1$ with probability $1/2$. Other interesting strategies for recovering a sparse signal in the time domain might be to sense a comparatively small number of its Fourier coefficients. In fact, [14] show that in this case our main condition holds with

$$S \asymp n/\log^6 p,$$

for nearly all subsets of observed coefficients of size $n$. Vershynin has informed us that $S \asymp n/\log^5 p$ also holds, and we believe that $S \asymp n/\log p$ is also true. More generally, suppose that $X$ is obtained by randomly sampling $n$ rows of a $p$ by $p$ orthonormal matrix $U$ (and renormalizing the columns so that they are unit-normed). Then we can take

$$S \asymp n/[\mu^2 \log^5 p],$$

with $\mu$ the coherence $\mu := \sup_{ij} \sqrt{n}|U_{ij}|$ [14].

Of course, all these calculations have implications for random designs. For example, suppose that in an *idealized application* one could—in a first experiment—observe $\beta$ directly and measure $y^{(1)} \sim N(\beta, \sigma^2 I_p)$. Consider a second experiment where one measures instead $y \sim N(X\beta, \sigma^2 I_n)$, where $X$ is a renormalized random design matrix with i.i.d. entries taking values $\pm 1$ with probability $1/2$. Suppose that the signal is $S$-sparse (note that $\|X\beta\|_{\ell_2} \asymp \|\beta\|_{\ell_2}$). Then reversing (2.1), we see that with about

$$n \asymp S \cdot \log(p/S)$$

observations, one would get just about the same mean squared error that one would achieve by measuring *all* the coordinates of $\beta$ directly (and applying thresholding).

Such procedures are not foreign to statisticians. Combining parameters by random design or otherwise goes back a long way; see, for example, the long history of blood pooling strategies. The theoretical analysis needs of course to be validated with numerical simulations, which may give further insights about the practical behavior of our methods. Section 4 presents a first series of experiments to complement our study.

2.2. *Applications.* The ability to recover a sparse or nearly sparse parameter vector from a few observations raises tantalizing opportunities and we mention just a few to give concrete ideas:

1. *Biomedical imaging.* In the field of biomedical imaging, one is often only able to collect far fewer measurements than the number of pixels. In magnetic resonance imaging (MRI), for instance, one would like to reconstruct high-resolution images from heavily undersampled frequency data, as this would allow image acquisition speeds far beyond those offered by current technologies; for example, see [33] and [16]. If the image is sparse, as is the case in magnetic resonance angiography (or if its gradient is sparse or, more generally, if the image is sparse in a fixed basis [9]), then $\ell_1$-minimization may have a chance to be very effective in such challenging settings.

2. *Analog to digital.* By making a number $n$ of general linear measurements rather than measuring the usual pixels, one could, in principle, reconstruct a compressible or sparse image with essentially the same resolution as one would obtain by measuring all the pixels. Now suppose one could design analog sensors able to make measurements by correlating the signal we wish to acquire against incoherent waveforms as discussed in the previous sections. Then one would effectively be able to make up a digital image with far fewer sensors than what is usually considered necessary [14, 19].

3. *Sensor networks.* There are promising applications in sensor networks where taking random projections may yield the same distortion (the same quality of reconstruction), but using much less power than what is traditionally required [29].

**3. Proof of theorems.** We now prove Theorems 1.1, 1.2 and 1.3, and we introduce some notation that we will use throughout this section. We let $X^1, \ldots, X^p \in \mathbf{R}^n$ be the $p$ columns of $X$ (the exploratory variables) so that $X\beta = \beta_1 X^1 + \cdots + \beta_p X^p$ and $(X^*y)_j = \langle y, X^j \rangle$, $1 \le j \le p$. We recall that the columns of $X$ are normalized to have unit norm, that is, $\|X^j\|_{\ell_2} = 1$.

Note that it is sufficient to prove our theorems with $\sigma = 1$, as the general case would follow from a simple rescaling argument. Therefore, we assume $\sigma = 1$ from now on. Now a key observation is that, with large probability, $z \sim N(0, I_n)$ obeys the *orthogonality condition*

$$(3.1) \qquad |\langle z, X^j \rangle| \le \lambda_p \qquad \text{for all } 1 \le j \le p,$$

for $\lambda_p = \sqrt{2 \log p}$. This is standard and simply follows from the fact that, for each $j$, $Z_j := \langle z, X^j \rangle \sim N(0,1)$. We will see that if (3.1) holds, then (1.10) holds. Note that, for each $u > 0$, $\mathbf{P}(\sup_j |Z_j| > u) \le 2p \cdot \phi(u)/u$, where $\phi(u) := (2\pi)^{-1/2} e^{-u^2/2}$, and our quantitative probabilistic statement just follows from this bound. Better bounds are possible, but we will not pursue these refinements here. As remarked earlier, if the columns were not unit normed, one would obtain the same conclusion with $\lambda_p = \sqrt{1 + \delta_1} \cdot \sqrt{2 \log p}$ since $\|X^j\|_{\ell_2} \le \sqrt{1 + \delta_1}$.

### 3.1. *High-dimensional geometry.*
It is probably best to start by introducing intuitive arguments underlying Theorems 1.1 and 1.2. These ideas are very geometrical and we hope they will convey the gist of the proof.

Consider Theorem 1.1 first, and suppose that $y = X\beta + z$, where $z$ obeys the orthogonality condition (3.1) for some $\lambda_p$. Let $\hat{\beta}$ be the minimizer of (1.7). Clearly, the true vector of parameters $\beta$ is feasible and, hence,

$$\|\hat{\beta}\|_{\ell_1} \le \|\beta\|_{\ell_1}.$$

Decompose $\hat{\beta}$ as $\hat{\beta} = \beta + h$ and let $T_0$ be the support of $\beta$, $T_0 = \{i : \beta_i \ne 0\}$. Then $h$ obeys two geometric constraints:

1. First, as essentially observed in [20],

$$\|\beta\|_{\ell_1} - \|h_{T_0}\|_{\ell_1} + \|h_{T_0^c}\|_{\ell_1} \le \|\beta + h\|_{\ell_1} \le \|\beta\|_{\ell_1},$$

where again the $i$th component of the vector $h_{T_0}$ is that of $h$ if $i \in T_0$ and zero otherwise (similarly for $h_{T_0^c}$). Hence, $h$ obeys the cone constraint

$$(3.2) \qquad \|h_{T_0^c}\|_{\ell_1} \le \|h_{T_0}\|_{\ell_1}.$$

2. Second, since

$$\langle z - r, X^j \rangle = \langle X\hat{\beta} - X\beta, X^j \rangle = \langle Xh, X^j \rangle,$$

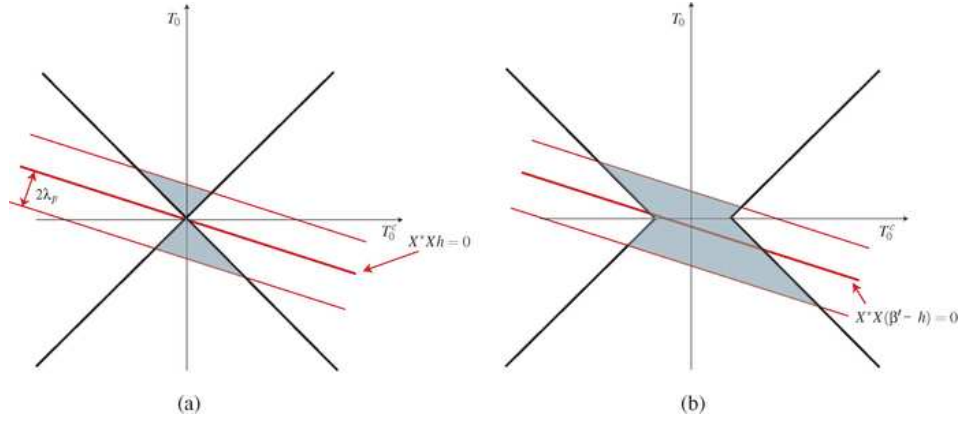

Fig. 1. *This figure represents the geometry of the constraints. On the left, the shaded area represents the set of $h$ obeying both* (3.2) *(hourglass region) and* (3.3) *(slab region). The right figure represents the situation in the more general case.*

it follows from the triangle inequality that

$$\|X^*Xh\|_{\ell_\infty} \leq 2\lambda_p. \tag{3.3}$$

We will see that these two geometrical constraints imply that $h$ is small in the $\ell_2$-norm. In other words, we will show that

$$\sup_{h \in \mathbf{R}^p} \|h\|_{\ell_2}^2 \quad \text{subject to} \quad \|h_{T_0^c}\|_{\ell_1} \leq \|h_{T_0}\|_{\ell_1} \quad \text{and} \quad \|X^*Xh\|_{\ell_\infty} \leq 2\lambda_p \tag{3.4}$$

obeys the desired bound, that is, $O(\lambda_p^2 \cdot |T_0|)$. This is illustrated in Figure 1(a). Hence, our statistical question is deeply connected with the geometry of high-dimensional Banach spaces, and that of high-dimensional spaces in general.

To think about the general case, consider the set of indices $T_0 := \{i : |\beta_i| > \sigma\}$ and let $\beta_{T_0}$ be the vector equal to $\beta$ on $T_0$ and zero outside, $\beta = \beta_{T_0} + \beta_{T_0^c}$. Suppose now that $\beta_{T_0}$ were feasible. Then we would have $\|\hat\beta\|_{\ell_1} \leq \|\beta_{T_0}\|_{\ell_1}$; writing $\hat\beta = \beta_{T_0} + h$, the same analysis as that of Theorem 1.1—and outlined above—would give

$$\|\hat\beta - \beta_{T_0}\|_{\ell_2}^2 = O(\log p) \cdot |T_0| \cdot \sigma^2.$$

From $\|\hat\beta - \beta\|_{\ell_2}^2 \leq 2\|\hat\beta - \beta_{T_0}\|_{\ell_2}^2 + 2\|\beta - \beta_{T_0}\|_{\ell_2}^2$, one would get

$$\|\hat\beta - \beta\|_{\ell_2}^2 = O(\log p) \cdot |T_0| \cdot \sigma^2 + 2 \sum_{i \,:\, |\beta_i| < \sigma} \beta_i^2,$$

which is the content of (1.13). Unfortunately, while $\beta_{T_0}$ may be feasible for "most" $S$-sparse vectors $\beta$, it is not for some, and the argument is considerably more involved.

### 3.2. *Proof of Theorem 1.1.*

LEMMA 3.1. *Suppose $T_0$ is a set of cardinality $S$ with $\delta + \theta < 1$. For a vector $h \in \mathbf{R}^p$, we let $T_1$ be the $S$ largest positions of $h$ outside of $T_0$. Put $T_{01} = T_0 \cup T_1$. Then*

$$\|h\|_{\ell_2(T_{01})} \le \frac{1}{1-\delta}\|X_{T_{01}}^T X h\|_{\ell_2} + \frac{\theta}{(1-\delta)S^{1/2}}\|h\|_{\ell_1(T_0^c)}$$

*and*

$$\|h\|_{\ell_2}^2 \le \|h\|_{\ell_2(T_{01})}^2 + S^{-1}\|h\|_{\ell_1(T_0^c)}^2.$$

PROOF. Consider the restricted transformation $X_{T_{01}} : \mathbf{R}^{T_{01}} \to \mathbf{R}^n$, $X_{T_{01}} c := \sum_{j \in T_{01}} c_j X^j$. Let $V \subset \mathbf{R}^n$ be the span of $\{X^j : j \in T_{01}\}$. Then $V$ is of course the range of $X_{T_{01}}$ and also the orthogonal complement of the kernel of $X_{T_{01}}^T$, which says that $\mathbf{R}^n$ is the orthogonal sum $V \oplus V^\perp$. Because $\delta < 1$, we know that the operator $X_{T_{01}}$ is a bijection from $\mathbf{R}^{T_{01}}$ to $V$, with singular values between $\sqrt{1-\delta}$ and $\sqrt{1+\delta}$. As a consequence, for any $c \in \ell_2(T_{01})$, we have

$$\sqrt{1-\delta}\|c\|_{\ell_2} \le \|X_{T_{01}} c\|_{\ell_2} \le \sqrt{1+\delta}\|c\|_{\ell_2}.$$

Moreover, letting $P_V$ denote the orthogonal projection onto $V$, we have for each $w \in \mathbf{R}^n$, $X_{T_{01}}^T w = X_{T_{01}}^T P_V w$ and it follows that

$$(3.5) \qquad \sqrt{1-\delta}\|P_V w\|_{\ell_2} \le \|X_{T_{01}}^T w\|_{\ell_2} \le \sqrt{1+\delta}\|P_V w\|_{\ell_2}.$$

We apply this to $w := Xh$ and conclude, in particular, that

$$(3.6) \qquad \|P_V X h\|_{\ell_2} \le (1-\delta)^{-1/2}\|X_{T_{01}}^T X h\|_{\ell_2}.$$

The next step is to derive a lower bound on $P_V X h$. To do this, we begin by dividing $T_0^c$ into subsets of size $S$ and enumerate $T_0^c$ as $n_1, n_2, \ldots, n_{p-|T_0|}$ in decreasing order of magnitude of $h_{T_0^c}$. Set $T_j = \{n_\ell, (j-1)S + 1 \le \ell \le jS\}$. That is, $T_1$ is as before and contains the indices of the $S$ largest coefficients of $h_{T_0^c}$, $T_2$ contains the indices of the next $S$ largest coefficients, and so on.

Decompose now $P_V X h$ as

$$(3.7) \qquad P_V X h = P_V X h_{T_{01}} + \sum_{j \ge 2} P_V X h_{T_j}.$$

By definition, $Xh_{T_{01}} \in V$ and $P_V X h_{T_{01}} = X h_{T_{01}}$. Further, since $P_V$ is an orthogonal projection onto the span of the $X^j$'s for $j \in T_{01}$, $P_V X h_{T_j} = \sum_{j \in T_{01}} c_j X^j$ for some coefficients $c_j$, and the following identity holds:

$$(3.8) \qquad \|P_V X h_{T_j}\|_{\ell_2}^2 = \langle P_V X h_{T_j}, X h_{T_j} \rangle.$$

By restricted orthogonality followed by restricted isometry, this gives

$$\langle P_V X h_{T_j}, X h_{T_j} \rangle \le \theta \left( \sum_{j \in T_{01}} |c_j|^2 \right)^{1/2} \|h_{T_j}\|_{\ell_2}$$

$$\le \frac{\theta}{\sqrt{1-\delta}} \|P_V X h_{T_j}\|_{\ell_2} \|h_{T_j}\|_{\ell_2},$$

which upon combining with (3.8) gives

$$(3.9) \qquad \|P_V X h_{T_j}\|_{\ell_2} \le \frac{\theta}{\sqrt{1-\delta}} \|h_{T_j}\|_{\ell_2}.$$

We then develop an upper bound on $\sum_{j \ge 2} \|h_{T_j}\|_{\ell_2}$ as in [12]. By construction, the magnitude of each component $h_{T_{j+1}}[i]$ of $h_{T_{j+1}}$ is less than the average of the magnitudes of the components of $h_{T_j}$,

$$|h_{T_{j+1}}[i]| \le \|h_{T_j}\|_{\ell_1}/S.$$

Then $\|h_{T_{j+1}}\|_{\ell_2}^2 \le \|h_{T_j}\|_{\ell_1}^2/S$ and, therefore,

$$(3.10) \qquad \sum_{j \ge 2} \|h_{T_j}\|_{\ell_2} \le S^{-1/2} \sum_{j \ge 1} \|h_{T_j}\|_{\ell_1} = S^{-1/2} \|h\|_{\ell_1(T_0^c)}.$$

To summarize, $X h_{T_{01}}$ obeys $\|X h_{T_{01}}\|_{\ell_2} \ge \sqrt{1-\delta} \|h_{T_{01}}\|_{\ell_2}$ by restricted isometry, and since $\sum_{j \ge 2} \|P_V X h_{T_j}\|_{\ell_2} \le \theta(1-\delta)^{-1/2} S^{-1/2} \|h\|_{\ell_1(T_0^c)}$,

$$\|P_V X h\|_{\ell_2} \ge \sqrt{1-\delta} \|h\|_{\ell_2(T_{01})} - \frac{\theta}{\sqrt{1-\delta}} S^{-1/2} \|h\|_{\ell_1(T_0^c)}.$$

Combining this with (3.6) proves the first part of the lemma.

It remains to argue the second part. Observe that the $k$th largest value of $h_{T_0^c}$ obeys

$$|h_{T_0^c}|_{(k)} \le \|h_{T_0^c}\|_{\ell_1}/k$$

and, therefore,

$$\|h_{T_{01}^c}\|_{\ell_2}^2 \le \|h_{T_0^c}\|_{\ell_1}^2 \sum_{k \ge S+1} 1/k^2 \le \|h_{T_0^c}\|_{\ell_1}^2/S,$$

which is what we needed to establish. The lemma is proven.  $\square$

Theorem 1.1 is now an easy consequence of this lemma. Observe that on the one hand, (3.2) gives

$$\|h_{T_0^c}\|_{\ell_1} \le \|h_{T_0}\|_{\ell_1} \le S^{1/2} \|h_{T_0}\|_{\ell_2},$$

while on the other hand, (3.3) gives

$$\|X_{T_{01}}^T X h\|_{\ell_2} \le (2S)^{1/2} \cdot 2\lambda_p,$$

since each of the $2S$ coefficients of $X_{T_{01}}^T X h$ is at most $2\lambda_p$ (3.3). In conclusion, we apply Lemma 3.1 and obtain

$$\|h\|_{\ell_2(T_{01})} \leq \frac{1}{1 - \delta - \theta} \cdot \sqrt{2S} \cdot 2\lambda_p.$$

The theorem follows since $\|h\|_{\ell_2}^2 \leq 2\|h\|_{\ell_2(T_{01})}^2$.

3.3. *Proof of Theorem* 1.3. The argument underlying Theorem 1.3 is almost the same as that of Theorem 1.1. We let $T_0$ be the set of the $S$ largest entries of $\beta$, and write $\hat{\beta} = \beta + h$ as before. If $z$ obeys the orthogonality condition (1.5), $\beta$ is feasible and

$$\|\beta_{T_0}\|_{\ell_1} - \|h_{T_0}\|_{\ell_1} + \|h_{T_0^c}\|_{\ell_1} - \|\beta_{T_0^c}\|_{\ell_1} \leq \|\beta + h\|_{\ell_1} \leq \|\beta\|_{\ell_1},$$

which gives

$$\|h_{T_0^c}\|_{\ell_1} \leq \|h_{T_0}\|_{\ell_1} + 2\|\beta_{T_0^c}\|_{\ell_1}.$$

The presence of the extra term is the only difference in the argument. We then conclude from Lemma 3.1 and (3.3) that

$$\|h_{T_0^c}\|_{\ell_2} \leq \frac{C}{1 - \delta - \theta} \cdot (\lambda_p \cdot S^{1/2} + \|\beta_{T_0^c}\|_{\ell_1} \cdot S^{-1/2}).$$

The second part of Lemma 3.1 gives $\|h\|_{\ell_2} \leq 2\|h\|_{\ell_2(T_{01})} + \|\beta_{T_0^c}\|_{\ell_1} \cdot S^{-1/2}$. Since for $\beta$ obeying the decay condition (1.16), $\|\beta_{T_0^c}\|_{\ell_1} \cdot S^{-1/2} \leq C \cdot R \cdot S^{-r}$, with $r = 1/s - 1/2$, we have established that, for all $S \leq S_*$,

$$\|h_{T_0^c}\|_{\ell_2} \leq \frac{C}{1 - \delta_{S_*} - \theta_{\theta_*, 2S_*}} \cdot (\lambda_p \cdot S^{1/2} + R \cdot S^{-r}).$$

The theorem follows.

3.4. *Proof of Theorem* 1.2. We begin with an auxiliary lemma.

LEMMA 3.2. *For any vector* $\beta$, *we have*

$$\|X\beta\|_{\ell_2} \leq \sqrt{1 + \delta}(\|\beta\|_{\ell_2} + (2S)^{-1/2}\|\beta\|_{\ell_1}).$$

PROOF. Let $T_1$ be the $2S$ largest positions of $\beta$, then $T_2$ be the next largest, and so forth. Then

$$\|X\beta\|_{\ell_2} \leq \|X\beta_{T_1}\|_{\ell_2} + \sum_{j \geq 2} \|X\beta_{T_j}\|_{\ell_2}.$$

From restricted isometry, we have

$$\|X\beta_{T_1}\|_{\ell_2} \leq (1 + \delta)^{1/2}\|\beta_{T_1}\|_{\ell_2} \leq (1 + \delta)^{1/2}\|\beta\|_{\ell_2}$$

and

$$\|X\beta_{T_j}\|_{\ell_2} \leq (1+\delta)^{1/2}\|\beta_{T_j}\|_{\ell_2} \leq (1+\delta)^{1/2}(2S)^{-1/2}\|\beta_{T_{j-1}}\|_{\ell_1}.$$

The claim follows. $\square$

We now turn to the proof of Theorem 1.2. As usual, we let $\hat{\beta}$ be the $\ell_1$ minimizer subject to the constraints

$$(3.11) \qquad \|X^*(X\hat{\beta} - y)\|_{\ell_\infty} = \sup_{1 \leq j \leq p} |\langle X\hat{\beta} - y, X^j\rangle| \leq (1 + t^{-1})\lambda,$$

where $\lambda := \sqrt{2\log p}$ for short.

Without loss of generality, we may order the $\beta_j$'s in decreasing order of magnitude

$$(3.12) \qquad |\beta_1| \geq |\beta_2| \geq \cdots \geq |\beta_p|.$$

In particular, by the sparsity assumption on $\beta$, we know that

$$(3.13) \qquad \beta_j = 0 \qquad \text{for all } j > S.$$

In particular, we see that

$$\sum_j \min(\beta_j^2, \lambda^2) \leq S \cdot \lambda^2.$$

Let $S_0$ be the smallest integer such that

$$(3.14) \qquad \sum_j \min(\beta_j^2, \lambda^2) \leq S_0 \cdot \lambda^2;$$

thus, $0 \leq S_0 \leq S$ and

$$(3.15) \qquad S_0 \cdot \lambda^2 \leq \lambda^2 + \sum_j \min(\beta_j^2, \lambda^2).$$

Also, observe from (3.12) that

$$S_0 \cdot \lambda^2 \geq \sum_{j=1}^{S_0+1} \min(\beta_j^2, \lambda^2) \geq (S_0 + 1)\min(\beta_{S_0+1}^2, \lambda^2)$$

and, hence, $\min(\beta_{S_0+1}^2, \lambda^2)$ is strictly less than $\lambda^2$. By (3.12), we conclude that

$$(3.16) \qquad \beta_j < \lambda \qquad \text{for all } j > S_0.$$

Write $\beta = \beta^{(1)} + \beta^{(2)}$, where

$$\beta_j^{(1)} = \beta_j \cdot 1_{1 \leq j \leq S_0},$$

$$\beta_j^{(2)} = \beta_j \cdot 1_{j > S_0}.$$

Thus, $\beta^{(1)}$ is a hard-thresholded version of $\beta$, localized to the set

$$T_0 := \{1, \ldots, S_0\}.$$

By (3.16), $\beta^{(2)}$ is $S$-sparse with

$$\|\beta^{(2)}\|_{\ell_2}^2 = \sum_{j > S_0} \min(\beta_j^2, \lambda^2) \leq S_0 \cdot \lambda^2.$$

As we shall see in the next section, Corollary A.3 allows the decomposition $\beta^{(2)} = \beta' + \beta''$, where

$$\|\beta'\|_{\ell_2} \leq \frac{1+\delta}{1-\delta-\theta} \lambda \cdot S_0^{1/2}, \tag{3.17}$$

$$\|\beta'\|_{\ell_1} \leq \frac{1+\delta}{1-\delta-\theta} \lambda \cdot S_0 \tag{3.18}$$

and

$$\|X^* X \beta''\|_{\ell_\infty} < \frac{1-\delta^2}{1-\delta-\theta} \lambda. \tag{3.19}$$

We use this decomposition and observe that

$$X^*(X(\beta^{(1)} + \beta') - y) = -X^* X \beta'' - X^* z$$

and, hence, by (3.1) and (3.19),

$$\|X^*(X(\beta^{(1)} + \beta') - y)\|_{\ell_\infty} \leq \left(1 + \frac{1-\delta^2}{1-\delta-\theta}\right)\lambda. \tag{3.20}$$

By assumption, $(1-\delta-\theta)^{-1} \leq t^{-1}$ and, therefore, $\beta^{(1)} + \beta'$ is feasible, which in turn implies

$$\|\hat{\beta}\|_{\ell_1} \leq \|\beta^{(1)} + \beta'\|_{\ell_1} \leq \|\beta^{(1)}\|_{\ell_1} + \frac{(1+\delta)}{1-\delta-\theta} S_0 \cdot \lambda.$$

Put $\hat{\beta} = \beta^{(1)} + h$. Then $\|\hat{\beta}\|_{\ell_1} \geq \|\beta^{(1)}\|_{\ell_1} - \|h\|_{\ell_1(T_0)} + \|h\|_{\ell_1(T_0^c)}$ so that

$$\|h\|_{\ell_1(T_0^c)} \leq \|h\|_{\ell_1(T_0)} + \frac{1+\delta}{1-\delta-\theta} S_0 \cdot \lambda, \tag{3.21}$$

and from (3.11) and (3.20), we conclude that

$$\|X^* X(\beta' - h)\|_{\ell_\infty} \leq 2\left(1 + \frac{1-\delta^2}{1-\delta-\theta}\right)\lambda. \tag{3.22}$$

Figure 1(b) schematically illustrates both these constraints.

The rest of the proof is essentially the same as that of Theorem 1.1. By Lemma 3.1, we have

$$\|h_{01}\|_{\ell_2} \leq \frac{1}{1-\delta}\|X_{T_{01}}^T Xh\|_{\ell_2} + \frac{\theta}{(1-\delta)S_0^{1/2}}\|h\|_{\ell_1(T_0^c)}.$$

On the other hand, from (3.22), we have

$$\|X_{T_{01}}^T X(\beta' - h)\|_{\ell_2} \leq 2\sqrt{2}\Big(1 + \frac{1-\delta^2}{1-\delta-\theta}\Big)S_0^{1/2} \cdot \lambda,$$

while from Lemma 3.2 and (3.18), (3.17), we have

$$\|X\beta'\|_{\ell_2} \leq (1+1/\sqrt{2})\frac{(1+\delta)^{3/2}}{1-\delta-\theta}S_0^{1/2} \cdot \lambda$$

and, hence, by restricted isometry,

$$\|X_{T_{01}}^T X\beta'\|_{\ell_2} \leq (1+1/\sqrt{2})\frac{(1+\delta)^2}{1-\delta-\theta}S_0^{1/2} \cdot \lambda.$$

In short,

$$\|X_{T_{01}}^T Xh\|_{\ell_2} \leq C_0 \cdot S_0^{1/2} \cdot \lambda,$$

where $C_0$ was defined in (1.14). We conclude that

$$\|h_{01}\|_{\ell_2} \leq \frac{C_0}{1-\delta}S_0^{1/2} \cdot \lambda + \frac{\theta}{(1-\delta)S_0^{1/2}}\|h\|_{\ell_1(T_0^c)}.$$

Finally, the bound (3.21) gives

$$\|h\|_{\ell_1(T_0^c)} \leq S_0^{1/2}\|h_{01}\|_{\ell_2} + \frac{1+\delta}{1-\delta-\theta}S_0 \cdot \lambda$$

and, hence,

$$\|h_{01}\|_{\ell_2} \leq C_0' \cdot S_0^{1/2} \cdot \lambda,$$

where

$$C_0' := \frac{C_0}{1-\delta-\theta} + \frac{\theta(1+\delta)}{(1-\delta-\theta)^2}.$$

Applying the second part of Lemma 3.1 and (3.21), we conclude

$$\|h\|_{\ell_2} \leq 2\|h_{01}\|_{\ell_2} + \frac{1+\delta}{1-\delta-\theta}S_0^{1/2} \cdot \lambda \leq C_2 \cdot S_0^{1/2} \cdot \lambda$$

and the claim follows from (3.15).

THE DANTZIG SELECTOR 23

3.5. *Refinements.* The constant $C_2$ obtained by this argument is not best possible; it is possible to lower it further, but at the cost of making the arguments slightly more complicated. For instance, in Lemma 3.2 one can exploit the approximate orthogonality between $X\beta_{T_1}$ and the $X\beta_{T_j}$'s to improve over the triangle inequality. Also, instead of defining $T_1, T_2, \ldots$ to have cardinality $S$, one can instead choose these sets to have cardinality $\rho S$ for some parameter $\rho$ to optimize in later. We will not pursue these refinements here. However, we observe that in the limiting case $\delta = \theta = 0$, then $X$ is an orthogonal matrix, and as we have seen earlier, $\hat{\beta}_j = \max(|(X^*y)_j| - \lambda, 0)\,\mathrm{sgn}((X^*y)_j)$. In this case one easily verifies that

$$\|\beta - \hat{\beta}\|_{\ell^2}^2 \leq \sum_{i=1}^{p} \min(\beta_i^2, 4\lambda^2).$$

This would correspond, roughly speaking, to a value of $C_2 = 2$ in Theorem 1.2, and therefore shows that there is room to improve $C_2$ by a factor of roughly 8.

## 4. Numerical experiments.

This section presents numerical experiments to illustrate the Dantzig selector and gives some insights about the numerical method for solving (1.9).

4.1. *An illustrative example.* In this first example, the design matrix $X$ has $n = 72$ rows and $p = 256$ columns, with independent Gaussian entries (and then normalized so that the columns have unit norm). We then select $\beta$ with $S := |\{i : \beta_i \neq 0\}| = 8$ and form $y = X\beta + z$, where the $z_i$'s are i.i.d. $N(0, \sigma^2)$. The noise level is adjusted so that

$$\sigma = \frac{1}{3}\sqrt{\frac{S}{n}}.$$

Here and below, the regularizing parameter $\lambda_p$ in $(DS)$ is chosen via Monte Carlo simulations, that is, as the empirical maximum of $|X^*z|_i$ over several realizations of $z \sim N(0, I_n)$. The results are presented in Figure 2.

First, we note that in this example our procedure correctly identifies all the nonzero components of $\beta$, and correctly sets to zero all the others. Quantitatively speaking, the ratio $\rho^2$ between the squared error loss and the ideal squared error (1.11) is equal to

(4.1) $$\rho^2 := \frac{\sum_i (\hat{\beta}_i - \beta_i)^2}{\sum_i \min(\beta_i^2, \sigma^2)} = 10.28.$$

(Note that here $2\log p = 11.09$.) Second, and as essentially observed earlier, the method clearly exhibits a soft-thresholding type of behavior and as a result, tends to underestimate the true value of the nonzero parameters. However, the two-stage Dantzig selector (1.20) introduced in Section 1.6 corrects for this bias. When applied to the same dataset, it recenters the estimator, and yields an improved squared error since now $\rho^2 = 1.14$; compare the results of Figures 2(a) and 2(b).

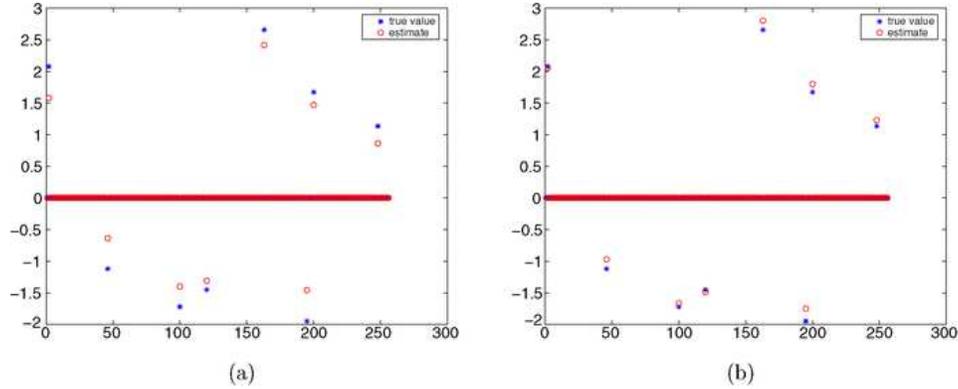

FIG. 2. *Estimation from $y = X\beta + z$ with $X$ a 72 by 256 matrix with independent Gaussian entries. A blue star indicates the true value of the parameter and a red circle the estimate. In this example, $\sigma = 0.11$ and $\lambda = 3.5$ so that the threshold is at $\delta = \lambda \cdot \sigma = 0.39$. (a) Dantzig selector* (1.7). *Note that our procedure correctly identifies all the nonzero components of $\beta$, and correctly sets to zero all the others. Observe the soft-thresholding-like behavior.* (b) *Estimation based on the two-stage strategy* (1.20). *The signal and estimator are very sparse, which is why there is a solid red line at zero.*

In our practical experience, the two-stage or Gauss–Dantzig selector procedure tends to outperform our original proposal, and to study its typical quantitative performance, we performed a series of experiments designed as follows:

1. $X$ is a 72 by 256 matrix, sampled as before ($X$ is fixed throughout);
2. select a support set $T$ of size $|T| = S$ uniformly at random, and sample a vector $\beta$ on $T$ with independent and identically distributed entries according to the model

$$\beta_i = \varepsilon_i(1 + |a_i|),$$

where the sign $\varepsilon_i = \pm 1$ with probability $1/2$, and $a_i \sim N(0,1)$ (the moduli and the signs of the amplitudes are independent);
3. make $\tilde{y} = X\beta + z$, with $z \sim N(0, \sigma^2 I_n)$ and compute $\hat{\beta}$ by means of the two-stage procedure (1.20);
4. repeat 500 times for each $S$, and for different noise levels $\sigma$.

The results are presented in Figure 3 and show that our approach works well. With the squared ratio $\rho^2$ as in (4.1), the median and the mean of $\rho^2$ are 2.35 and 9.42, respectively, for a noise level $\sigma$ set at $1/3\sqrt{S/n}$. In addition, 75% of the time $\rho^2$ is less than 10. With $\sigma = \sqrt{S/n}$, the mean and median are 12.38 and 13.78, respectively, with a bell-shaped distribution.

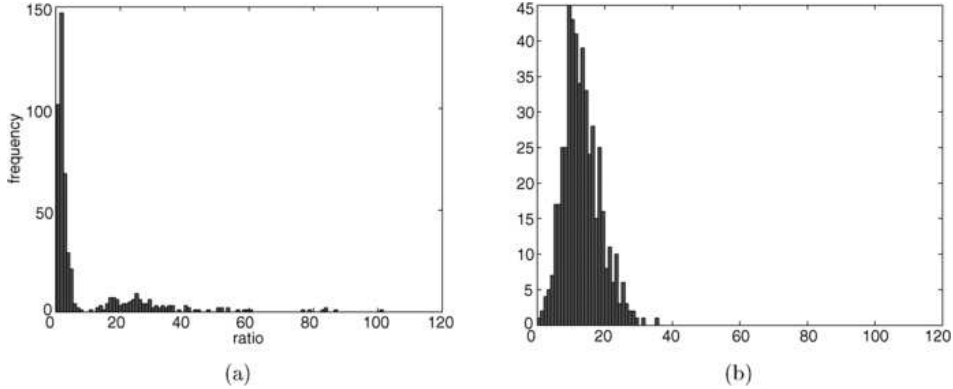

Fig. 3. *Statistics of the ratio between the squared error $\sum_i (\hat{\beta}_i - \beta_i)^2$ and the ideal mean-squared error $\sum_i \min(\beta_i^2, \sigma^2)$. (a) $S = 8$, $\sigma = 1/3\sqrt{S/n} = 0.11$, and the threshold is here $\lambda \cdot \sigma = 0.5814$. (b) $S = 8$, $\sigma = \sqrt{S/n} = 0.33$, and the threshold is now $\lambda \cdot \sigma = 1.73$.*

The reader may wonder why the two histograms in Figure 3 look somewhat different. The answer is simply due to the fact that in Figure 3(a) $\sigma = 0.11$ and the threshold is about $\lambda \cdot \sigma = 0.5814$, which means that the nonzero $\beta_i$'s are above the noise level. In Figure 3(b), however, $\sigma = 0.33$, and the threshold is $\lambda \cdot \sigma = 1.73$. This means that a fraction of the nonzero components of $\beta$ are within the noise level and will be set to zero. This explains the observed difference in the quantitative behavior.

4.2. *Binary design matrices.* We now consider the case where the design matrix has i.i.d. entries taking on values $\pm 1$, each with probability $1/2$ (the entries are then divided by $\sqrt{n}$ so that the columns have unit norm). This simply amounts to measuring differences among randomly selected subsets of coordinates. Note that if $X$ had $0/1$ entries, one would measure the aggregate effect of randomly selected subsets of coordinates, much like in pooling design schemes. The number of predictors is set to $p = 5,000$ and the number of observations to $n = 1,000$. Of interest here is the estimation accuracy as the number $S$ of significant parameters increases. The results are presented in Table 1. In all these experiments, the nonzero coordinates of the parameter vector are sampled as in Section 4.1 and the noise level is adjusted to $\sigma = 1/3\sqrt{S/n}$, so that with $y = X\beta + z$, the variance $\mathbf{E}\|z\|^2 = n\sigma^2$ is proportional to $\mathbf{E}\|\beta\|^2$ with the same constant of proportionality (fixed signal-to-noise ratio).

*Ratio between the squared error $\sum_i(\hat{\beta}_i - \beta_i)^2$ and the ideal mean squared error $\sum_i \min(\beta_i^2, \sigma^2)$. The binary matrix $X$ is the same in all these experiments, and the noise level $\sigma = 1/3\sqrt{S/n}$ is adjusted so that the signal-to-noise ratio is nearly constant. Both estimators and especially the Gauss–Dantzig selector exhibit a remarkable performance until a breakdown point around $S = 200$*

| $S$ | 5 | 10 | 20 | 50 | 100 | 150 | 200 |
|---|---|---|---|---|---|---|---|
| Dantzig selector | 22.81 | 17.30 | 28.85 | 18.49 | 25.71 | 49.73 | 74.93 |
| Gauss–Dantzig selector | 0.36 | 0.65 | 1.04 | 1.09 | 1.53 | 13.71 | 48.74 |

Our estimation procedures and most notably the Gauss–Dantzig selector are remarkably accurate as long as the number of significant parameters is not too large, here about $S = 200$. For example, for $S = 100$, the ratio between the Gauss–Dantzig selector's squared error loss and the ideal mean-squared error is only 1.53. Figure 4 illustrates the estimation precision in this case.

To confirm these findings, we now sample the amplitudes of the parameter vector $\beta$ according to a Cauchy distribution in order to have a wide range of component values $\beta_i$, some of which are within the noise level, while others are way above; $X$ is fixed and we now vary the number $S$ of nonzero components of $\beta$ as before, while $\sigma = 0.5$ is now held constant. The results are presented in Table 2. Again, the Gauss–Dantzig selector performs well.

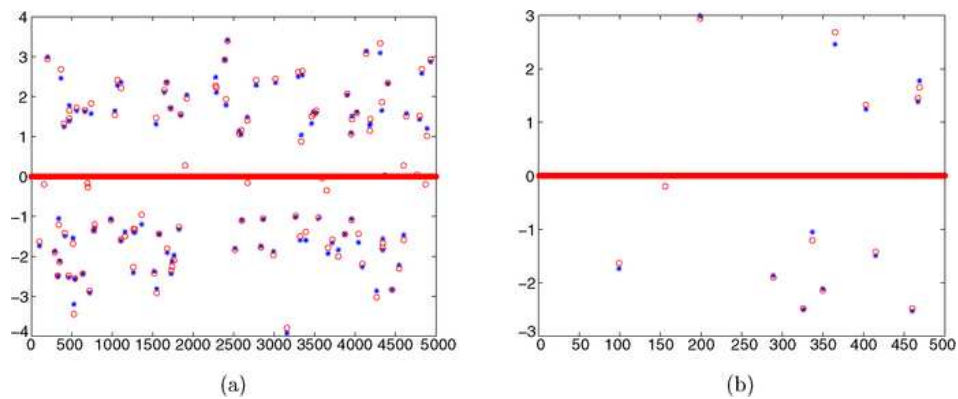

Fig. 4. *Estimation from $y = X\beta + z$ with $X$ a 1,000 by 5,000 matrix with independent binary entries. A blue star indicates the true value of the parameter and a red circle the estimate.* (a) *True parameter values and estimates obtained via the Dantzig selector (1.20). There are $S = 100$ significant parameters.* (b) *Same plot but showing the first 500 coordinates for higher visibility.*

TABLE 2
*Ratio between the squared error $\sum_i (\hat{\beta}_i - \beta_i)^2$ and the ideal mean squared error $\sum_i \min(\beta_i^2, \sigma^2)$. The binary matrix $X$, $\sigma = 0.5$ and $\lambda \cdot \sigma = 2.09$ are the same in all these experiments*

| $S$ | 5 | 10 | 20 | 50 | 100 | 150 | 200 |
|---|---|---|---|---|---|---|---|
| Gauss–Dantzig selector | 3.70 | 4.52 | 2.78 | 3.52 | 4.09 | 6.56 | 5.11 |

4.3. *Examples in signal processing.* We are interested in recovering a one-dimensional signal $f \in \mathbf{R}^p$ from noisy and undersampled Fourier coefficients of the form

$$y_j = \langle f, \phi_j \rangle + z_j, \qquad 1 \leq j \leq n,$$

where $\phi_j(t)$, $t = 0, \ldots, p-1$, is a sinusoidal waveform $\phi_j(t) = \sqrt{2/n} \cos(\pi(k_j + 1/2)(t + 1/2))$, $k_j \in \{0, 1, \ldots, p-1\}$. Consider the signal $f$ in Figure 5; $f$ is not sparse, but its wavelet coefficients sequence $\beta$ is. Consequently, we may just as well estimate its coefficients in a nice wavelet basis. Letting $\Phi$ be the matrix with the $\phi_k$'s as rows, and $W$ be the orthogonal wavelet matrix with wavelets as columns, we have $y = X\beta + z$, where $X = \Phi W$, and our estimation procedure applies as is.

The test signal is of size $p = 4{,}096$ (Figure 5), and we sample a set of frequencies of size $n = 512$ by extracting the lowest 128 frequencies and randomly selecting the others. With this set of observations, the goal is to study the quantitative behavior of the Gauss–Dantzig selector procedure for various noise levels. [Owing to the factor $\sqrt{2/n}$ in the definition of $\phi_j(t)$, the columns of $X$ have size about one and for each column, individual thresholds

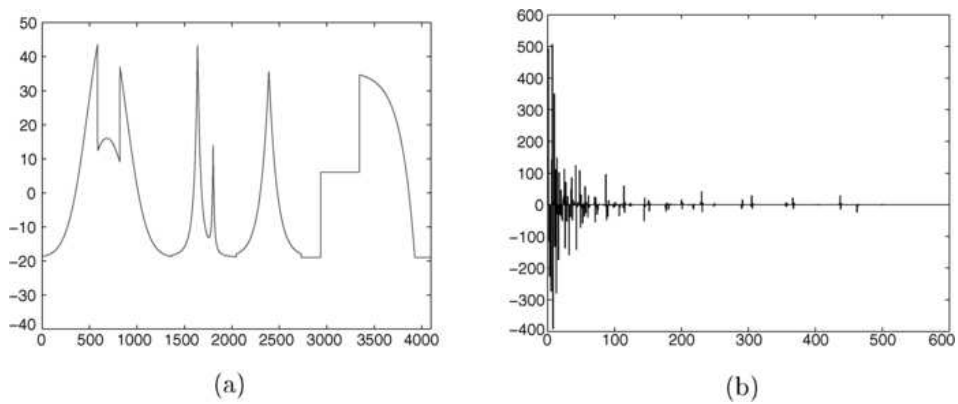

FIG. 5.   (a) *One-dimensional signal $f$ we wish to reconstruct.* (b) *First* 512 *wavelet coefficients of $f$.*

TABLE 3

*Performance of the Gauss–Dantzig procedure in estimating a signal from undersampled and noisy Fourier coefficients. The subset of variables is here estimated by $|\hat{\beta}_i| > \sigma/4$, with $\hat{\beta}$ as in* (1.7). *The top row is the value of the signal-to-noise ratio (SNR)*

| SNR $\alpha = \|X\beta\|/\sqrt{n\sigma^2}$ | 100 | 20 | 10 | 2 | 1 | 0.5 |
|---|---|---|---|---|---|---|
| $\sum_i(\hat{\beta}_i - \beta_i)^2 / \sum_i \min(\beta_i^2, \sigma^2)$ | 15.51 | 2.08 | 1.40 | 1.47 | 0.91 | 1.00 |

$\lambda_i$—$|(X^*r)_i| \leq \lambda_i \cdot \sigma$—are determined by looking at the empirical distribution of $|(X^*z)_i|$.] We adjust $\sigma$ so that $\alpha^2 = \|X\beta\|_{\ell_2}^2 / \mathbf{E}\|z\|_{\ell_2}^2 = \|X\beta\|_{\ell_2}^2 / n\sigma^2$ for various levels of the signal-to-noise ratio $\alpha$. We use Daubechies' wavelets with four vanishing moments for the reconstruction. The results are presented in Table 3. As one can see, high statistical accuracy holds over a wide range of noise levels. Interestingly, the estimator is less accurate when the noise level is very small ($\alpha = 100$), which is not surprising, since in this case there are 178 wavelet coefficients exceeding $\sigma$ in absolute value.

In our last example, we consider the problem of reconstructing an image from undersampled Fourier coefficients. Here $\beta(t_1, t_2)$, $0 \leq t_1, t_2 < N$, is an unknown $N$ by $N$ image so that $p$ is the number of unknown pixels, $p = N^2$. As usual, the data is given by $y = X\beta + z$, where

$$(4.2) \quad (X\beta)_k = \sum_{t_1,t_2} \beta(t_1, t_2) \cos(2\pi(k_1 t_1 + k_2 t_2)/N), \qquad k = (k_1, k_2),$$

or $(X\beta)_k = \sum_{t_1,t_2} \beta(t_1, t_2) \sin(2\pi(k_1 t_1 + k_2 t_2)/N)$. In our example [see Figure 6(b)], the image $\beta$ is not sparse, but the gradient is. Therefore, to reconstruct the image, we apply our estimation strategy and minimize the $\ell_1$-norm of the gradient size, also known as the total-variation of $\beta$,

$$(4.3) \qquad\qquad \min \|\tilde{\beta}\|_{\text{TV}} \quad \text{subject to} \quad |(X^*r)_i| \leq \lambda_i \cdot \sigma$$

(the individual thresholds again depend on the column sizes as before); formally, the total-variation norm is of the form

$$\|\tilde{\beta}\|_{\text{TV}} = \sum_{t_1,t_2} \sqrt{|D_1\tilde{\beta}(t_1, t_2)|^2 + |D_2\tilde{\beta}(t_1, t_2)|^2},$$

where $D_1$ is the finite difference $D_1\tilde{\beta} = \beta(t_1, t_2) - \beta(t_1 - 1, t_2)$ and $D_2\tilde{\beta} = \tilde{\beta}(t_1, t_2) - \tilde{\beta}(t_1, t_2 - 1)$; in short, $\|\tilde{\beta}\|_{\text{BV}}$ is the $\ell_1$-norm of the size of the gradient $D\tilde{\beta} = (D_1\tilde{\beta}, D_2\tilde{\beta})$; see also [34].

Our example follows the data acquisition patterns of many real imaging devices which can collect high-resolution samples along radial lines at relatively few angles. Figure 6(a) illustrates a typical case where one gathers $N = 256$ samples along each of 22 radial lines. In a first experiment then, we observe $22 \times 256$ noisy real-valued Fourier coefficients and use (4.3) for the recovery problem illustrated in Figure 6. The number of observations is then $n = 5{,}632$, whereas there are $p = 65{,}536$ observations. In other words, about 91.5% of the 2D Fourier coefficients of $\beta$ are missing. The SNR in this experiment is equal to $\|X\beta\|_{\ell_2}/\|z\|_{\ell_2} = 5.85$. Figure 6(c) shows the reconstruction obtained by setting the unobserved Fourier coefficients to zero, while (d) shows the reconstruction (4.3). We follow up with a second experiment where the unknown image is now 512 by 512 so that $p = 262{,}144$ and $n = 22 \times 512 = 11{,}264$. The fraction of missing Fourier coefficients is now approaching 96%. The SNR ratio is about the same, $\|X\beta\|_{\ell_2}/\|z\|_{\ell_2} = 5.77$.

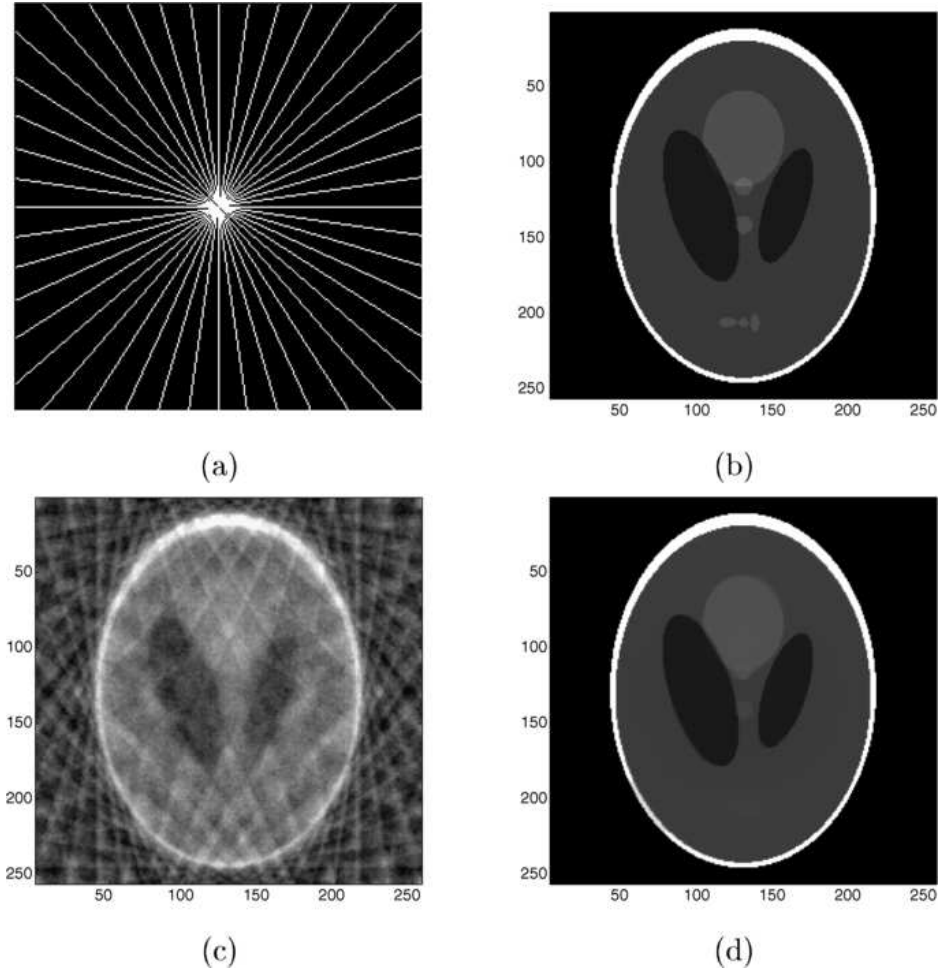

Fig. 6.  (a) *Sampling "domain" in the frequency plane; Fourier coefficients are sampled along* 22 *approximately radial lines; here,* $n \approx 0.086p$. (b) *The Logan–Shepp phantom test image.* (c) *Minimum energy reconstruction obtained by setting unobserved Fourier coefficients to zero.* (d) *Reconstruction obtained by minimizing the total-variation, as in* (4.3).

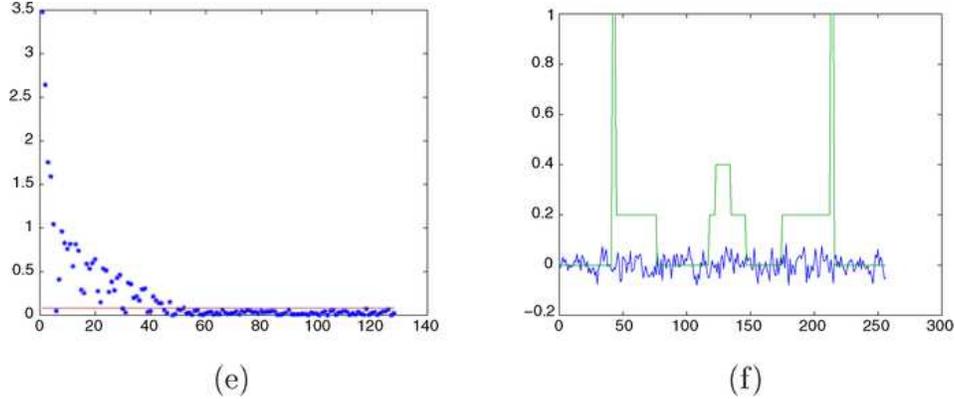

(e)                                      (f)

FIG. 6. *Continued.* (e) *Magnitude of the true Fourier coefficients along a radial line (frequency increases from left to right) on a logarithmic scale. Blue stars indicate values of* $\log(1 + |(X\beta)_k|)$, *while the solid red line indicates the noise level* $\log(1+\sigma)$. *Less than a third of the frequency samples exceed the noise level.* (f) $X^*z$ *and* $\beta$ *are plotted along a scanline to convey a sense of the noise level.*

The reconstructions are of very good quality, especially when compared to the naive reconstruction which minimizes the energy of the reconstruction subject to matching the observed data. Figure 7 also shows the middle horizontal scanline of the phantom. As expected, we note a slight loss of contrast due to the nature of the estimator which here operates by "soft-thresholding" the gradient. There are, of course, ways of correcting for this bias, but such issues are beyond the scope of this paper.

4.4. *Implementation.* In all the experiments above, we used a primal-dual interior point algorithm for solving the linear program (1.9). We used a specific implementation which we outline, as this gives some insight about the computational workload of our method. For a general linear program with inequality constraints

$$\min c^*\beta \quad \text{subject to} \quad F\beta \le b,$$

define

- $f(\beta) = X\beta - b$,
- $r_{\text{dual}} = c + X^*\lambda$,
- $r_{\text{cent}} = -\text{diag}(\lambda)f(\beta) - \mathbf{1}/t$,

where $\lambda \in \mathbf{R}^m$ are the so-called dual variables, and $t$ is a parameter whose value typically increases geometrically at each iteration; there are as many dual variables as inequality constraints. In a standard primal-dual method (with logarithmic barrier function) [8], one updates the current primal-dual

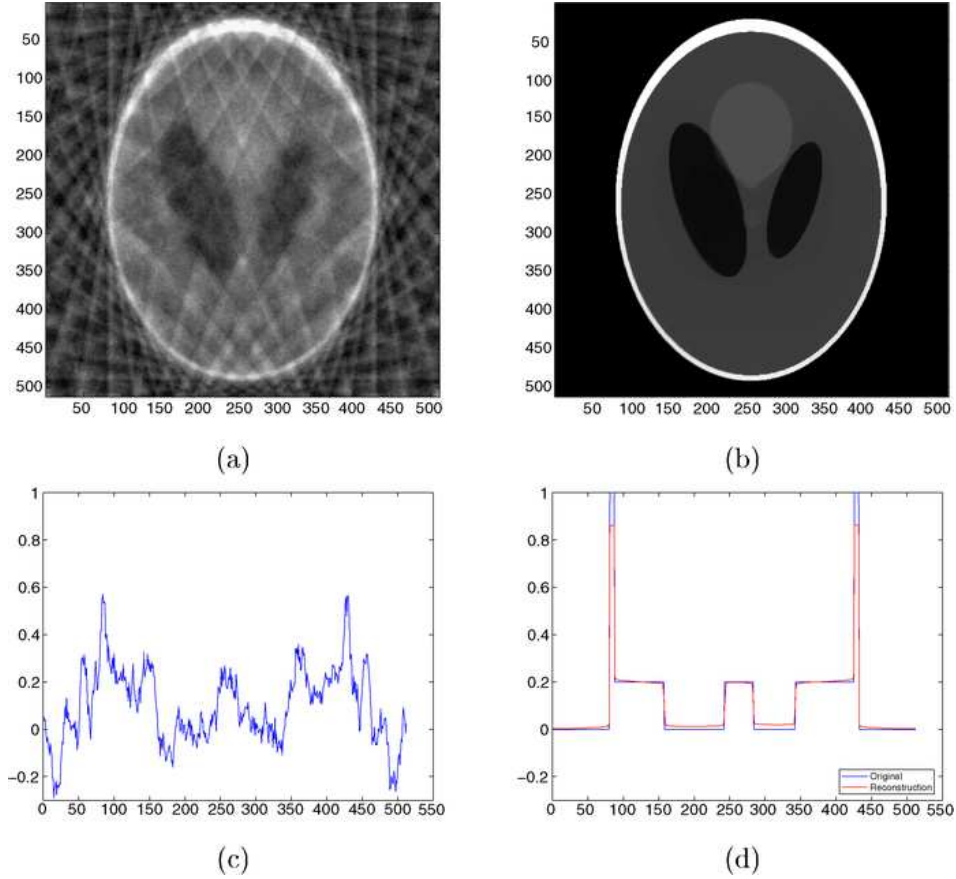

Fig. 7. (a) and (b). *Similar experience as in Figure* 6 *but at a higher resolution* $(p = 512^2)$ *so that now* $n \approx 0.043p$. (c) *and* (d). *Scanlines of both reconstructions.*

pair $(\beta, \lambda)$ by means of a Newton step, and solves

$$\begin{pmatrix} 0 & F^* \\ -\operatorname{diag}(\lambda)F & -\operatorname{diag}(f(\beta)) \end{pmatrix} \begin{pmatrix} \Delta\beta \\ \Delta\lambda \end{pmatrix} = -\begin{pmatrix} r_{\text{dual}} \\ r_{\text{cent}} \end{pmatrix};$$

that is,

$$\Delta\lambda = -\operatorname{diag}(1/f(\beta))(\operatorname{diag}(\lambda)F\Delta_\beta - r_{\text{cent}}),$$

$$-[F^* \operatorname{diag}(\lambda/f(\beta))F]\Delta\beta = -(r_{\text{dual}} + F^* \operatorname{diag}(1/f(\beta))r_{\text{cent}}).$$

The current guess is then updated via $(\beta^+, \lambda^+) = (\beta, \lambda) + s(\Delta\beta, \Delta s)$, where the stepsize $s$ is determined by line search or otherwise. Typically, the sequence of iterations stops once the primal-dual gap and the size of the residual vector fall below a specified tolerance level [8].

Letting $U = X^*X$ and $\tilde{y} = X^*y$ in (1.9), our problem parameters have the block structure

$$F = \begin{pmatrix} I & -I \\ -I & -I \\ U & 0 \\ -U & 0 \end{pmatrix}, \qquad b = \begin{pmatrix} 0 \\ 0 \\ \delta + \tilde{y} \\ \delta - \tilde{y} \end{pmatrix}, \qquad c = \begin{pmatrix} 0 \\ 1 \end{pmatrix},$$

which gives

$$F^* \operatorname{diag}(\lambda/f)F = \begin{pmatrix} D_1 + D_2 + U^*(D_3 + D_4)U & D_2 - D_1 \\ D_2 - D_1 & D_1 + D_2 \end{pmatrix},$$

where $D_i = \operatorname{diag}(\lambda_i/f_i)$, $1 \le i \le 4$, and

$$\begin{pmatrix} f_1 \\ f_2 \\ f_3 \\ f_4 \end{pmatrix} = \begin{pmatrix} \tilde{\beta} - u \\ -\tilde{\beta} - u \\ U\tilde{\beta} - \delta - \tilde{y} \\ -U\tilde{\beta} - \delta + \tilde{y} \end{pmatrix}.$$

Put $\binom{r_1}{r_2} = r_{\text{dual}} + F^* \operatorname{diag}(1/f(\beta))r_{\text{cent}}$. It is convenient to solve the system

$$F^* \operatorname{diag}(\lambda/f)F \begin{pmatrix} \Delta\tilde{\beta} \\ \Delta u \end{pmatrix} = \begin{pmatrix} r_1 \\ r_2 \end{pmatrix}$$

by elimination and obtain

$$(4(D_1 + D_2)^{-1}D_1 D_2 + U^*(D_3 + D_4)U)\Delta\tilde{\beta} = r_1 - (D_1 + D_2)^{-1}(D_1 - D_2)r_2,$$

$$(D_1 + D_2)\Delta u = r_2 - (D_2 - D_1)\Delta\tilde{\beta}.$$

In other words, each step may involve solving a $p$ by $p$ system of linear equations. In fact, when $n$ is less than $p$, it is possible to solve a smaller system thanks to the Sherman–Woodbury–Morrison formula. Indeed, write $U^*(D_3 + D_4)U = X^*B$, where $B = X(D_3 + D_4)X^*X$ and put $D_{12} = 4(D_1 + D_2)^{-1}D_1 D_2$. Then

$$(D_{12} + X^T B)^{-1} = D_{12}^{-1} - D_{12}^{-1}X^*(I + BD_{12}^{-1}X^*)^{-1}BD_{12}^{-1}.$$

The advantage is that one needs to solve the smaller $n$ by $n$ system $(I + BD_{12}^{-1}X^*)\beta' = b'$. Hence, the cost of each Newton iteration is essentially that of solving an $n$ by $n$ system of linear equations, plus that of forming the matrix $(I + BD_{12}^{-1}X^*)$. As far as the number of Newton iterations is concerned, we ran thousands of experiments and have never needed more than 45 Newton iterations for convergence.

Note that in some important applications, we may have fast algorithms for applying $X$ and $X^*$ to an arbitrary vector, as in the situation where $X$ is a partial Fourier matrix, since one can make use of FFT's. In such settings, one never forms $X^*X$ and uses iterative algorithms such as Conjugate Gradients for solving such linear systems; of course, this speeds up the computations.

Finally, a collection of MATLAB routines solving (1.9) for reproducing some of these experiments and testing these ideas in other setups is available at the address www.l1-magic.org/.

## 5. Discussion.

5.1. *Significance for model selection.* We would like to briefly discuss the implications of this work for model selection. Given a data matrix $X$ (with unit-normed columns) and observations $y = X\beta + z$, our procedure will estimate $\beta$ by that vector with minimum $\ell_1$-norm among all objects $\tilde{\beta}$ obeying

(5.1) $$|X^*r|_i \leq (1 + t^{-1})\sqrt{2\log p} \cdot \sigma, \qquad \text{where } r = y - X\tilde{\beta}.$$

As the theory and the numerical experiments suggest, many of the coordinates $\hat{\beta}_i$ will typically be zero (at least under the assumption that the true unknown vector $\beta$ is sparse) and, therefore, our estimation procedure effectively returns a candidate model $\hat{I} := \{i : \hat{\beta}_i \neq 0\}$.

As we have seen, the Dantzig selector is guaranteed to produce optimal results if

(5.2) $$\delta(X)_{2S} + \theta(X)_{S,2S} < 1$$

[note that since $\theta(X)_{S,2S} \leq \delta(X)_{3S}$, it would be sufficient to have $\delta(X)_{2S} + \delta(X)_{3S} < 1$]. We have commented on the interpretation of this condition already. In a typical model selection problem, $X$ is given; it is then natural to ask if this particular $X$ obeys (5.2) for the assumed level $S$ of sparsity. Unfortunately, obtaining an answer to this question might be computationally prohibitive, as it may require checking the extremal singular values of exponentially many submatrices.

While this may represent a limitation, two observations are in order. First, there are empirical evidence and theoretical analysis suggesting approximate answers for certain types of random matrices, and there is nowadays a significant amount of activity developing tools to address these questions [17]. Second, the failure of (5.2) to hold is in general indicative of a structural difficulty of the problem, so that any procedure is likely to be unsuccessful for sparsity levels in the range of the critical one. To illustrate this point, let us consider the following example. Suppose that $\delta_{2S} + \theta_{S,2S} > 1$. Then this says that $\delta_{2S}$ is large and since $\delta_S$ is increasing in $S$, it may very well be that for $S'$ in the range of $S$, for example, $S' = 3S$, there might be submatrices (in fact, possibly many submatrices) with $S'$ columns which are either rank-deficient or which have very small singular values. In other words, the significant entries of the parameter vector might be arranged in such a way so that even if one knew their location, one would not be able to estimate their values because of rank-deficiency. This informally suggests that the Dantzig selector breaks down near the point where any estimation procedure, no matter how intractable, would fail.

Finally, it is worth emphasizing the connections between RIC [26] and (1.7). Both methods suggest a penalization which is proportional to the logarithm of the number of explanatory variables—a penalty that is well justified theoretically [7]. For example, in the very special case where $X$ is orthonormal, $p = n$, RIC applies a hard-thresholding rule to the vector $X^*y$ at about $O(\sqrt{2 \log p})$, while our procedure translates in a soft-thresholding at about the same level; in our convex formulation (1.7), this threshold level is required to make sure that the true vector is feasible. In addition, the ideas developed in this paper have broad applicability, and it is likely that they might be deployed and give similar bounds for RIC variable selection. Despite such possible similarities, our method differs substantially from RIC in terms of computational effort since (1.7) is tractable while RIC is not. In fact, we are not aware of any work in the model selection literature which is close in intent and in the results.

5.2. *Connections with other work.* In [10] a related problem is studied where the goal is to recover a vector $\beta$ from incomplete and contaminated measurements $y = X\beta + e$, where the error vector (stochastic or deterministic) obeys $\|e\|_{\ell_2}^2 \leq D$. There, the proposed reconstruction searches, among all objects consistent with the data $y$, for that with minimum $\ell_1$-norm,

$$(5.3) \qquad (P_2) \qquad \min \|\tilde{\beta}\|_{\ell_1} \quad \text{subject to} \quad \|y - X\tilde{\beta}\|_{\ell_2}^2 \leq D.$$

Under essentially the same conditions as those of Theorems 1.1 and 1.2, [10] showed that the reconstruction error is bounded by

$$(5.4) \qquad \qquad \|\beta^\sharp - \beta\|_{\ell_2}^2 \ \leq \ C_3^2 \cdot D.$$

In addition, it is also argued that, for arbitrary errors, the bound (5.4) is sharp and that, in general, one cannot hope for a better accuracy.

If one assumes that the error is stochastic as in this paper, however, a mean squared error of about $D$ might be far from optimal. Indeed, with the linear model (1.1), $D \sim \sigma^2 \chi_n$ and, therefore, $D$ has size about $n\sigma^2$. But suppose now $\beta$ is sparse and has only three nonzero coordinates, say, all exceeding the noise level. Then whereas Theorem 1.1 gives a loss of about $3\sigma^2$ (up to a log factor), (5.4) only guarantees an error of size about $n\sigma^2$. What is missing in [10] and is achieved here is the adaptivity to the unknown level of sparsity of the object we try to recover. Note that we do not claim that the program $(P_2)$ is ill-suited for adaptivity. It is possible that refined arguments would yield estimators based on quadratically constrained $\ell_1$-minimization [variations of (5.3)] obeying the special adaptivity properties discussed in this paper.

Last but not least, and while working on this manuscript, we became aware of related work [29]. Motivated by recent results [11, 14, 18], the authors studied the problem of reconstructing a signal from noisy random projections and obtained powerful quantitative estimates which resemble ours.

Their setup is different though, since they exclusively work with random design matrices, in fact, random Rademacher projections, and do not study the case of fixed $X$'s. In contrast, our model for $X$ is deterministic and does not involve any kind of randomization, although our results can of course be specialized to random matrices. But more importantly, and perhaps this is the main difference, their estimation procedure requires solving a combinatorial problem much like (1.15), whereas we use linear programming.

## APPENDIX

This appendix justifies the construction of a pseudo-hard thresholded vector which obeys the constraints; see (3.17), (3.18) and (3.19) in the proof of Theorem 1.2.

LEMMA A.1 (Dual sparse reconstruction, $\ell_2$ version).    *Let $S$ s.t. $\delta + \theta < 1$, and let $c_T$ be supported on $T$ for some $|T| \leq 2S$. Then there exist $\beta$ supported on $T$, and an exceptional set $E$ disjoint from $T$ with*

$$(6.1) \qquad |E| \leq S,$$

*such that*

$$(6.2) \qquad \langle X\beta, X^j \rangle = c_j \qquad \text{for all } j \in T$$

*and*

$$(6.3) \qquad |\langle X\beta, X^j \rangle| \leq \frac{\theta}{(1-\delta)\sqrt{S}} \|c_T\|_{\ell_2} \qquad \text{for all } j \notin (T \cup E)$$

*and*

$$(6.4) \qquad \left( \sum_{j \in E} |\langle X\beta, X^j \rangle|^2 \right)^{1/2} \leq \frac{\theta}{1-\delta} \|c_T\|_{\ell_2}.$$

*Also we have*

$$(6.5) \qquad \|\beta\|_{\ell_2} \leq \frac{1}{1-\delta} \|c_T\|_{\ell_2}$$

*and*

$$(6.6) \qquad \|\beta\|_{\ell_1} \leq \frac{\sqrt{2S}}{1-\delta} \|c_T\|_{\ell_2}.$$

PROOF.    We define $\beta$ by

$$\beta_T := (X_T^T X_T)^{-1} c_T,$$

and zero outside of $T$, which gives (6.2), (6.5) and (6.6) by Cauchy–Schwarz. Note that $X\beta = X_T\beta_T$. We then set

$$E := \left\{ j \notin T : |\langle X\beta, X^j\rangle| > \frac{\theta}{(1-\delta)\sqrt{S}}\|c_T\|_{\ell_2} \right\},$$

so (6.3) holds.

Now if $T'$ is disjoint from $T$ with $|T'| \leq S$ and $d_{T'}$ is supported on $T'$, then

$$|\langle X_{T'}^T X_T\beta_T, d_{T'}\rangle| \leq \theta\|\beta_T\|_{\ell_2}\|d_{T'}\|_{\ell_2}$$

and, hence, by duality,

$$(6.7) \qquad \|X^* X_T\beta_T\|_{\ell_2(T')} \leq \theta\|\beta_T\|_{\ell_2} \leq \frac{\theta}{1-\delta}\|c_T\|_{\ell_2}.$$

If $|E| \geq S$, then we can find $T' \subset E$ with $|T'| = S$. But then we have

$$\|X^* X_T\beta_T\|_{\ell_2(T')} > \frac{\theta}{(1-\delta)S^{1/2}}\|c_T\|_{\ell_2}|T'|^{1/2},$$

which contradicts (6.7). Thus, we have (6.1). Now we can apply (6.7) with $T' := E$ to obtain (6.4).  □

COROLLARY A.2 (Dual sparse reconstruction, $\ell_\infty$ version).  *Let $c_T$ be supported on $T$ for some $|T| \leq S$. Then there exists $\beta$ obeying (6.2) such that*

$$(6.8) \qquad |\langle X\beta, X^j\rangle| \leq \frac{\theta}{(1-\delta-\theta)\sqrt{S}}\|c_T\|_{\ell_2} \qquad \text{for all } j \notin T.$$

*Furthermore, we have*

$$(6.9) \qquad \|\beta\|_{\ell_2} \leq \frac{1}{1-\delta-\theta}\|c_T\|_{\ell_2}$$

*and*

$$(6.10) \qquad \|\beta\|_{\ell_1} \leq \frac{\sqrt{2S}}{1-\delta-\theta}\|c_T\|_{\ell_2}.$$

PROOF.  The proof of this lemma operates by iterating the preceding lemma as in Lemma 2.2 of [13]. We simply rehearse the main ingredients and refer the reader to [13] for details.

We may normalize $\sum_{j\in T} |c_j|^2 = 1$. Write $T_0 := T$ and note that $|T_0| \leq 2S$. Using Lemma A.1, we can find a vector $\beta^{(1)}$ and a set $T_1 \subseteq \{1, \ldots, p\}$ such

that

$$T_0 \cap T_1 = \varnothing,$$

$$|T_1| \leq S,$$

$$\langle X\beta^{(1)}, X^j \rangle = c_j \qquad \text{for all } j \in T_0,$$

$$|\langle X\beta^{(1)}, X^j \rangle| \leq \frac{\theta}{(1-\delta)S^{1/2}} \qquad \text{for all } j \notin T_0 \cup T_1,$$

$$\left( \sum_{j \in T_1} |\langle X\beta^{(1)}, X^j \rangle|^2 \right)^{1/2} \leq \frac{\theta}{1-\delta},$$

$$\|\beta^{(1)}\|_{\ell_2} \leq \frac{1}{1-\delta},$$

$$\|\beta^{(1)}\|_{\ell_1} \leq \frac{\sqrt{2S}}{1-\delta}.$$

Applying Lemma A.1 iteratively gives a sequence of vectors $\beta^{(n+1)} \in \mathbf{R}^p$ and sets $T_{n+1} \subseteq \{1, \ldots, p\}$ for all $n \geq 1$ with the properties

$$T_{n+1} \cap (T_0 \cup T_n) = \varnothing,$$

$$|T_{n+1}| \leq S,$$

$$\langle X\beta^{(n+1)}, X^j \rangle = \langle X\beta^{(n)}, X^j \rangle \qquad \text{for all } j \in T_n,$$

$$\langle X\beta^{(n+1)}, X^j \rangle = 0 \qquad \text{for all } j \in T_0,$$

$$|\langle X\beta^{(n+1)}, X^j \rangle| \leq \frac{\theta}{(1-\delta)S^{1/2}} \left( \frac{\theta}{1-\delta} \right)^n \qquad \forall j \notin T_0 \cup T_n \cup T_{n+1},$$

$$\left( \sum_{j \in T_{n+1}} |\langle X\beta^{(n+1)}, X^j \rangle|^2 \right)^{1/2} \leq \frac{\theta}{1-\delta} \left( \frac{\theta}{1-\delta} \right)^n,$$

$$\|\beta^{(n+1)}\|_{\ell_2} \leq \frac{1}{1-\delta} \left( \frac{\theta}{1-\delta} \right)^n,$$

$$\|\beta^{(n+1)}\|_{\ell_1} \leq \frac{\sqrt{2S}}{1-\delta} \left( \frac{\theta}{1-\delta} \right)^n.$$

By hypothesis, we have $\frac{\theta}{1-\delta} \leq 1$. Thus, if we set

$$\beta := \sum_{n=1}^{\infty} (-1)^{n-1} \beta^{(n)},$$

then the series is absolutely convergent and, therefore, $\beta$ is a well-defined vector. And it turns out that $\beta$ obeys the desired properties; see Lemma 2.2 in [13].  $\square$

COROLLARY A.3 (Constrained thresholding). *Let $\beta$ be $S$-sparse such that*

$$\|\beta\|_{\ell_2} < \lambda \cdot S^{1/2}$$

*for some $\lambda > 0$. Then there exists a decomposition $\beta = \beta' + \beta''$ such that*

$$\|\beta'\|_{\ell_2} \leq \frac{1+\delta}{1-\delta-\theta}\|\beta\|_{\ell_2},$$

$$\|\beta'\|_{\ell_1} \leq \frac{1+\delta}{1-\delta-\theta}\frac{\|\beta\|_{\ell_2}^2}{\lambda}$$

*and*

$$\|X^* X \beta''\|_{\ell_\infty} < \frac{1-\delta^2}{1-\delta-\theta}\lambda.$$

PROOF. Let

$$T := \{j : |\langle X\beta, X^j \rangle| \geq (1+\delta)\lambda\}.$$

Suppose that $|T| \geq S$. Then we can find a subset $T'$ of $T$ with $|T'| = S$. Then by restricted isometry we have

$$(1+\delta)^2\lambda^2 S \leq \sum_{j \in T'} |\langle X\beta, X^j \rangle|^2 \leq (1+\delta)\|X\beta\|_{\ell_2}^2 \leq (1+\delta)^2\|\beta\|_{\ell_2}^2,$$

contradicting the hypothesis. Thus, $|T| < S$. Applying restricted isometry again, we conclude

$$(1+\delta)^2\lambda^2 |T| \leq \sum_{j \in T} |\langle X\beta, X^j \rangle|^2 \leq (1+\delta)^2\|\beta\|_{\ell_2}^2$$

and, hence,

$$|T| \leq S := \frac{\|\beta\|_{\ell_2}^2}{\lambda^2}.$$

Applying Corollary A.2 with $c_j := \langle X\beta, X^j \rangle$, we can find a $\beta'$ such that

$$\langle X\beta', X^j \rangle = \langle X\beta, X^j \rangle \qquad \text{for all } j \in T,$$

$$\|\beta'\|_{\ell_2} \leq \frac{1+\delta}{1-\delta-\theta}\|\beta\|_{\ell_2},$$

$$\|\beta'\|_{\ell_1} \leq \frac{(1+\delta)\sqrt{S}}{1-\delta-\theta}\|\beta\|_{\ell_2} = \frac{(1+\delta)}{1-\delta-\theta}\frac{\|\beta\|_{\ell_2}^2}{\lambda}$$

*and*

$$|\langle X\beta', X^j \rangle| \leq \frac{\theta(1+\delta)}{(1-\delta-\theta)\sqrt{S}}\|\beta\|_{\ell_2} \qquad \text{for all } j \notin T.$$

By definition of $S$, we thus have

$$|\langle X\beta', X^j\rangle| \le \frac{\theta}{1-\delta-\theta}(1+\delta)\lambda \qquad \text{for all } j \notin T.$$

Meanwhile, by definition of $T$, we have

$$|\langle X\beta, X^j\rangle| < (1+\delta)\lambda \qquad \text{for all } j \notin T.$$

Setting $\beta'' := \beta - \beta'$, the claims follow. $\quad\square$

**Acknowledgments.** Emmanuel Candès thanks Rob Nowak for sending him an early preprint, Hannes Helgason for bibliographical research on this project, Justin Romberg for his help with numerical simulations and Anestis Antoniadis for comments on an early version of the manuscript. We also thank the referees for their helpful remarks.

## REFERENCES

[1] Akaike, H. (1974). A new look at the statistical model identification. *IEEE Trans. Automatic Control* **19** 716–723. MR0423716

[2] Antoniadis, A. and Fan, J. (2001). Regularization of wavelet approximations (with discussion). *J. Amer. Statist. Assoc.* **96** 939–967. MR1946364

[3] Baraud, Y. (2000). Model selection for regression on a fixed design. *Probab. Theory Related Fields* **117** 467–493. MR1777129

[4] Barron, A. R., Birgé, L. and Massart, P. (1999). Risk bounds for model selection via penalization. *Probab. Theory Related Fields* **113** 301–413. MR1679028

[5] Barron, A. R. and Cover, T. M. (1991). Minimum complexity density estimation. *IEEE Trans. Inform. Theory* **37** 1034–1054. MR1111806

[6] Birgé, L. and Massart, P. (1997). From model selection to adaptive estimation. In *Festschrift for Lucien Le Cam* (D. Pollard, E. Torgersen and G. L. Yang, eds.) 55–87. Springer, New York. MR1462939

[7] Birgé, L. and Massart, P. (2001). Gaussian model selection. *J. Eur. Math. Soc.* **3** 203–268. MR1848946

[8] Boyd, S. and Vandenberghe L. (2004). *Convex Optimization.* Cambridge Univ. Press. MR2061575

[9] Candès, E. J. and Romberg, J. (2005). Practical signal recovery from random projections. In *Computational Imaging III: Proc. SPIE International Symposium on Electronic Imaging* **1** 76–86. San Jose, CA.

[10] Candès, E. J., Romberg, J. and Tao, T. (2006). Stable signal recovery from incomplete and inaccurate measurements. *Comm. Pure Appl. Math.* **59** 1207–1223. MR2230846

[11] Candès, E. J., Romberg, J. and Tao, T. (2006). Robust uncertainty principles: Exact signal reconstruction from highly incomplete frequency information. *IEEE Trans. Inform. Theory* **52** 489–509. MR2236170

[12] Candès, E. J., Rudelson, M., Vershynin, R. and Tao, T. (2005). Error correction via linear programming. In *Proc. 46th Annual IEEE Symposium on Foundations of Computer Science (FOCS)* 295–308. IEEE, Los Alamitos, CA.

[13] Candès, E. J. and Tao, T. (2005). Decoding by linear programming. *IEEE Trans. Inform. Theory* **51** 4203–4215. MR2243152

[14] CANDÈS, E. J. and TAO, T. (2006). Near-optimal signal recovery from random projections: Universal encoding strategies? *IEEE Trans. Inform. Theory* **52** 5406–5425. MR2300700

[15] CHEN, S. S., DONOHO, D. L. and SAUNDERS, M. A. (1998). Atomic decomposition by basis pursuit. *SIAM J. Sci. Comput.* **20** 33–61. MR1639094

[16] DANIEL, B. L., YEN, Y. F., GLOVER, G. H. et al. (1998). Breast disease: Dynamic spiral MR imaging. *Radiology* **209** 499–509.

[17] DAUBECHIES, I. (2005). Personal communication.

[18] DONOHO, D. L. (2006). For most large underdetermined systems of linear equations the minimal $\ell_1$-norm solution is also the sparsest solution. *Comm. Pure Appl. Math.* **59** 797–829. MR2217606

[19] DONOHO, D. L. (2006). Compressed sensing. *IEEE Trans. Inform. Theory* **52** 1289–1306. MR2241189

[20] DONOHO, D. L. and HUO, X. (2001). Uncertainty principles and ideal atomic decomposition. *IEEE Trans. Inform. Theory* **47** 2845–2862. MR1872845

[21] DONOHO, D. L. and JOHNSTONE, I. M. (1994). Ideal spatial adaptation by wavelet shrinkage. *Biometrika* **81** 425–455. MR1311089

[22] DONOHO, D. L. and JOHNSTONE, I. M. (1994). Ideal denoising in an orthonormal basis chosen from a library of bases. *C. R. Acad. Sci. Paris Sér. I Math.* **319** 1317–1322. MR1310679

[23] DONOHO, D. L. and JOHNSTONE, I. M. (1995). Empirical atomic decomposition. Unpublished manuscript.

[24] ELAD, M. and BRUCKSTEIN, A. M. (2002). A generalized uncertainty principle and sparse representation in pairs of bases. *IEEE Trans. Inform. Theory* **48** 2558–2567. MR1929464

[25] FAN, J. and PENG, H. (2004). Nonconcave penalized likelihood with a diverging number of parameters. *Ann. Statist.* **32** 928–961. MR2065194

[26] FOSTER, D. P. and GEORGE, E. I. (1994). The risk inflation criterion for multiple regression. *Ann. Statist.* **22** 1947–1975. MR1329177

[27] FUCHS, J. (2004). On sparse representations in arbitrary redundant bases. *IEEE Trans. Inform. Theory* **50** 1341–1344. MR2094894

[28] GREENSHTEIN, E. and RITOV, Y. (2004). Persistence in high-dimensional linear predictor selection and the virtue of overparametrization. *Bernoulli* **10** 971–988. MR2108039

[29] HAUPT, J. and NOWAK, R. (2006). Signal reconstruction from noisy random projections. *IEEE Trans. Inform. Theory* **52** 4036–4048. MR2298532

[30] KETTENRING, J., LINDSAY, B. and SIEGMUND, D., eds. (2003). Statistics: Challenges and opportunities for the twenty-first century. NSF report. Available at www.pnl.gov/scales/docs/nsf_report.pdf.

[31] MALLOWS, C. L. (1973). Some comments on $C_P$. *Technometrics* **15** 661–675.

[32] NATARAJAN, B. K. (1995). Sparse approximate solutions to linear systems. *SIAM J. Comput.* **24** 227–234. MR1320206

[33] PETERS, D. C., KOROSEC, F. R., GRIST, T. M., BLOCK, W. F., HOLDEN, J. E., VIGEN, K. K. and MISTRETTA, C. A. (2000). Undersampled projection reconstruction applied to MR angiography. *Magnetic Resonance in Medicine* **43** 91–101.

[34] RUDIN, L. I., OSHER, S. and FATEMI, E. (1992). Nonlinear total variation based noise removal algorithms. *Physica D* **60** 259–268.

[35] Sardy, S., Bruce, A. G. and Tseng, P. (2000). Block coordinate relaxation methods for nonparametric wavelet denoising. *J. Comput. Graph. Statist.* **9** 361–379. MR1822091

[36] Schwarz, G. (1978). Estimating the dimension of a model. *Ann. Statist.* **6** 461–464. MR0468014

[37] Szarek, S. J. (1991). Condition numbers of random matrices. *J. Complexity* **7** 131–149. MR1108773

[38] Tibshirani, R. (1996). Regression shrinkage and selection via the lasso. *J. Roy. Statist. Soc. Ser. B* **58** 267–288. MR1379242

Applied and Computational Mathematics, MC 217–50
California Institute of Technology
Pasadena, California 91125
USA
E-mail: emmanuel@acm.caltech.edu

Department of Mathematics
University of California
Los Angeles, California 90095
USA
E-mail: tao@math.ucla.edu